\documentclass[11pt,letterpaper]{amsart}

\usepackage[dvipsnames,svgnames,x11names,hyperref]{xcolor}
\usepackage[hyphens]{url}
\usepackage[colorlinks=true,linkcolor=Maroon,citecolor=blue,urlcolor=blue,hypertexnames=false,linktocpage]{hyperref}
\usepackage{bookmark}
\usepackage{amsmath,thmtools,mathtools}
\mathtoolsset{showonlyrefs=true}

\usepackage{fancyhdr}
\usepackage{esint}
\usepackage{enumerate}

\usepackage{pictexwd,dcpic}
\usepackage{graphicx}
\usepackage{caption}

\newcounter{mgncount}

\swapnumbers
\declaretheorem[name=Theorem,numberwithin=section]{thm}
\declaretheorem[name=Remark,style=remark,sibling=thm]{rem}
\declaretheorem[name=Lemma,sibling=thm]{lemma}
\declaretheorem[name=Proposition,sibling=thm]{prop}
\declaretheorem[name=Conditions,sibling=thm]{cond}

\declaretheorem[name=Corollary,sibling=thm]{cor}

\numberwithin{equation}{section}

\newcommand{\ti}{\tilde}

\newcommand{\cn}{\colon}
\newcommand{\sub}{\subset}
\newcommand{\ov}{\overline}
\newcommand{\mr}{\mathring}
\newcommand{\rr}{\mathtt r}

\newcommand{\bbN}{\mathbb{N}}

\newcommand{\bbR}{\mathbb{R}}

\newcommand{\bbS}{\mathbb{S}}

\newcommand{\bbM}{\mathbb{M}}

\newcommand{\8}{\infty}

\newcommand{\al}{\alpha}
\newcommand{\be}{\beta}
\newcommand{\ga}{\gamma}
\newcommand{\de}{\delta}
\newcommand{\ep}{\epsilon}
\newcommand{\ka}{\kappa}
\newcommand{\la}{\lambda}
\newcommand{\om}{\omega}
\newcommand{\si}{\sigma}
\newcommand{\Si}{\Sigma}

\newcommand{\vp}{\varphi}
\newcommand{\ze}{\zeta}

\newcommand{\Om}{\Omega}
\newcommand{\De}{\Delta}
\newcommand{\Ga}{\Gamma}

\newcommand{\cH}{\mathcal{H}}

\newcommand{\cO}{\mathcal{O}}

\newcommand{\cS}{\mathcal{S}}

\newcommand{\cU}{\mathcal{U}}

\newcommand{\del}{\partial}
\newcommand{\n}{\nabla}

\newcommand{\fa}{\forall}
\newcommand{\rt}{\sqrt}

\newcommand{\co}{{\rm c}_{_K}}
\newcommand{\s}{{\rm s}_{_K}}

\newcommand{\ip}[2]{\left\langle #1,#2 \right\rangle}
\newcommand{\fr}[2]{\frac{#1}{#2}}
\newcommand{\tfr}[2]{\tfrac{#1}{#2}}
\newcommand{\x}{\times}

\DeclareMathOperator{\id}{id}

\DeclareMathOperator{\dist}{dist}

\DeclareMathOperator{\sgn}{sgn}

\newcommand{\pf}[1]{\begin{proof}#1 \end{proof}}
\newcommand{\eq}[1]{\begin{equation}\begin{alignedat}{2} #1 \end{alignedat}\end{equation}}

\newcommand{\br}[1]{\left(#1\right)}
\newcommand{\abs}[1]{\lvert #1\rvert}
\newcommand{\enum}[1]{\begin{enumerate}[(i)] #1 \end{enumerate}}

\newcommand{\ra}{\rightarrow}
\newcommand{\hra}{\hookrightarrow}

\newcommand{\mt}{\mapsto}

\newcommand{\mc}{\mathcal}

\newcommand{\mrm}{\mathrm}

\newcommand{\hp}{\hphantom}
\newcommand{\q}{\quad}

\begin{document}
\title[Quermassintegral preserving mean curvature flow]{The quermassintegral preserving mean curvature flow in the sphere }
\author[E. Cabezas-Rivas]{Esther Cabezas-Rivas}
\address{Universitat de València, Departament de Matemàtiques, Dr Moliner.~50, 46100 Burjassot, València (Spain)}
\email{Esther.Cabezas-Rivas@uv.es}
\author[J. Scheuer]{Julian Scheuer}
\address{Goethe-Universit\"at, Institut f\"ur Mathematik, Robert-Mayer-Str.~10, 60325 Frankfurt, Germany}
\email{scheuer@math.uni-frankfurt.de}
\date{\today. The work of the first author is partially supported by the AEI (Spain) and FEDER project PID2019-105019GB-C21, and by the GVA project AICO 2021 21/378.01/1. }
\begin{abstract}
We introduce a mean curvature flow with global term of convex hypersurfaces in the sphere, for which the global term can be chosen to keep any quermassintegral fixed.  Then, starting from a strictly convex initial hypersurface, we prove that the flow exists for all times and converges smoothly to a geodesic sphere. This provides a workaround to an issue present in the volume preserving mean curvature flow in the sphere introduced by Huisken in 1987.  We also classify solutions for some constant curvature type equations in space forms, as well as solitons in the sphere and in the upper branch of the De Sitter space.
\end{abstract}

\subjclass[2020]{53E10, 53C21.}

\maketitle

\section{Introduction and statement of Main Results}

Let $n \geq 2$ and $M^n \sub \bbM_{K}^{n+1}$ be a smooth, closed, embedded  hypersurface in a simply connected space form $\bbM_{K}^{n+1}$ of constant curvature $K \in \mathbb R$, given by the embedding $x_{0}$. We consider a family of embeddings $x = x(t,\cdot)$ satisfying the following mean curvature type flow with a global forcing term 
\eq{\label{flow-eq}\partial_t{x} = (\mu(t){\rm c}_{_{K}}(r) - H)\nu,}
which has initial condition $x(0, \cdot)=x_0$. Here $H$ is the mean curvature and $\nu$ the outward unit normal of the evolving hypersurfaces $M_t$.  For convex hypersurfaces (i.e. with $\kappa_1 \geq 0$, where $\kappa_1 \leq \cdots \leq \kappa_n$ denote the principal curvatures), the sign conventions are taken so that $-H \nu$ points inwards. Moreover, let $r$ denote the radial distance to a given point $\cO\in \bbM^{n+1}_{K}$, which we call {\it{origin}} in the sequel. This means that the flow \eqref{flow-eq} depends on the choice of the origin and, in fact, along the flow we will change the origin in a controlled way. We use the  notation 
\eq{\label{trig} {\rm c}_{_K}(r)=  {\rm s}_{_K}'(r),\quad \text{where} \quad {\rm s}_{_K}(r)= { \begin{cases} K^{-\frac1{2}}\sin(\sqrt{K}
	r),  & \text{ if } K>0 \\
r, &  \text{ if }  K =0 \\
 |K|^{-\frac1{2}} \sinh(\sqrt{|K|}
	r), &  \text{ if } K < 0.
\end{cases}}}

If $\sigma_\ell$ represents the $\ell$-th elementary symmetric function, we define the time-dependent term by
\eq{\label{def-mu} \mu(t) = \frac{\int_M  H \sigma_\ell \, d V_t}{\int_M  {\rm c}_{_{K}} \sigma_\ell \, d V_t},}
for each $\ell = 0, 1, \ldots, n$, where $d V_t$ denotes the volume element of $M_t$.
This choice guarantees that \eqref{flow-eq} gives a family of globally constrained mean curvature flows, where $\mu$ can be chosen to preserve any of the $n+1$ quermassintegrals $W_\ell(\Omega_t)$ of the evolving hypersurfaces $M_{t}$ (see \Autoref{sec:quermass} for a review of the quermassintegrals). Here $\Om_{t}$ denotes the convex region enclosed by $M_{t} = x(t,\bbS^{n})$ (note that we may assume that the common domain of the embeddings is $\bbS^{n}$, due to convexity).

Let us stress that for $\ell = 0$ and $K = 0$ the flow \eqref{flow-eq} with non-local term as in \eqref{def-mu} coincides with the {\it volume-preserving mean curvature flow} ({\sc vpmcf}) introduced by Huisken in 1987. He proved that strictly convex hypersurfaces in $\mathbb R^{n+1}$ remain convex and embedded under the flow, the solution exists for all times and converges to a round sphere smoothly as $t \to \infty$. Since then it was still an open question of extending the result to an $(n+1)$-dimensional sphere $\bbS^{n+1}_{K}$, $K>0$, where convexity can be lost under {\sc vpmcf}, as pointed out by Huisken \cite[page~38]{Huisken:/1987}.

Our main result settles this question by proposing the flow \eqref{flow-eq} as the most natural generalization of the {\sc vpmcf} to a space form with positive curvature. Indeed, such a definition preserves convexity under the flow, and allows us to prove the following version of Huisken's original result within the half sphere, where for a point $p\in \bbS^{n+1}_{K}$, $\cH(p)$ denotes the open hemisphere around $p$. 

\begin{thm}\label{thm:main}
Let $n\geq 2$ and $M_{0}\sub\bbS^{n+1}_{K}$ be a strictly convex hypersurface enclosing a domain $\Om_{0}$. Then there exists a finite system of origins $(\cO_{i})_{0\leq i\leq m}$  and numbers $0=t_{0}<t_{1}<\dots<t_{m}<t_{m+1}=\8$, such that the problem
\eq{\label{thm:main-A}\del_{t}x &= (\mu_{i}(t){\rm c}_{_{K}}(r_{i})-H)\nu,\q t\in [t_{i},t_{i+1}),\q 0\leq i\leq m,\\
		x(0,M) &= M_{0}\\
		x(t_{i},M)& = \lim_{t\nearrow t_{i}}M_{t},\q 1\leq i\leq m, }
where $r_{i}$ is the distance to $\cO_{i}$ and $\mu_{i}$ is given as in \eqref{def-mu} to keep the quermassintegral $W_\ell(\Omega_t)$ for any $\ell=0, 1, \ldots, n$ fixed, has a solution 
\eq{x\cn [0,\8)\x \bbS^{n}\ra \bbS_{K}^{n+1}.} For every $t\geq 0$, the embeddings $x(t,\cdot)$ smoothly map $\bbS^{n}$ to strictly convex hypersurfaces with
\eq{\cO_{i}\in \Om_{t}\q\mbox{and}\q M_{t}\sub \cH(\cO_{i})\q\fa t\in[t_{i},t_{i+1})} and satisfy spatial $C^{\8}$-estimates which are uniform in time. The restriction
\eq{x\cn [t_{m},\8)\x \bbS^{n}\ra \bbS_{K}^{n+1}}
is smooth and converges for $t\ra \8$ in $C^{\8}$ to a geodesic sphere around $\cO_{m}$ with radius determined by $W_{\ell}(B_{r}) = W_{\ell}(\Om_{0})$.
\end{thm}

At this stage we should mention that Guan and Li \cite{GuanLi:/2015} invented a {\it purely local} mean curvature type flow in the sphere, which is volume-preserving and drives starshaped hypersurfaces to geodesic spheres.  There is also a flow of Guan-Li type that preserves  $W_{\ell-1}(\Omega_t)$ and decreases $W_{\ell}(\Omega_t)$, which has so far refused to allow curvature estimates (see \cite{ChenGuanLiScheuer:/2022} for an overview over known results). 

However, notice that non-local flows are more challenging than their corresponding local counterparts, as the evolution depends heavily on the global shape of the hypersurface $M_t$ and the presence of the term $\mu(t)$ in all the relevant evolution equations causes a plethora of extra complications; e.g. comparison principles and preservation of key properties fail (cf.~\cite{CaMi4}), as well as embedded hypersurfaces may develop self-intersections (cf.~\cite{MayerSimonett:/2000}).

In this framework, our flow \eqref{flow-eq} is, to the best of our knowledge, the first known  curvature flow in the sphere, which preserves any desired quermassintegral by a suitable choice of the global term $\mu$ and which enjoys smooth convergence to a geodesic sphere.

\subsection*{Elliptic counterpart: rigidity results}

As a byproduct of the pinching estimates (see \autoref{ev-p-parabolic}) that we require to prove \Autoref{thm:main}, we can also classify hypersurfaces in a space form $\bbM_{K}^{n+1}$   which have a rotationally symmetric curvature function, under suitable assumptions on the sectional curvature.

With this goal, we work with curvature functions more general than $H$. Let $\Ga\sub\bbR^{n}$ be a symmetric, open cone containing the positive cone 
\eq{\Ga_{+} = \{\ka\in \bbR^{n}\cn \ka_{i}>0, \text{ for all } i=1, \ldots,n\},}
and consider a symmetric function $f \in C^2(\Gamma)$. Let 
\eq{F(A) = f(\ka_{1},\dots,\ka_{n})}
be the corresponding operator dependent function, where $A$ denotes the Weingarten or shape operator. We assume the following:
\begin{cond} \label{cond} Let $f(\kappa) = F(A)$ be a $C^2$ symmetric function, defined on an open, symmetric cone $\Gamma \supset \Gamma_+$. We ask further that
	\begin{enumerate}
		\item[{\rm (a)}] $f$ is strictly increasing in each argument.
		\item[{\rm (b)}] $f$ is homogeneous of degree 1.
		\item[{\rm (c)}] $f$ is normalized so that $f(1, \cdots, 1) = n$. 
		
	\end{enumerate} 
\end{cond}
Notice that (a) implies that $F$ defines a strictly elliptic operator on $M$, as proved in \cite{HuiskenPolden:/1996}. We say that $f$ is {\it inverse concave/convex} if the dual function 
\[\tilde f(\ka_{1},\dots,\ka_{n}) = f(\ka_{1}^{-1},\dots,\ka_{n}^{-1})\]
is concave/convex  (see \autoref{symmetric} for a more detailed introduction).
 
A classical result by Alexandrov \cite{Alexandroff:12/1962} says that if a compact hypersuface embedded in $\mathbb R^{n+1}$ has $H =$ constant, it
must be a round sphere. Later on, Ros \cite{Ros:/1987} extended this result to the constancy of higher order symmetric functions $\sigma_\ell$. Hypersurfaces in a model space $\bbM_{K}^{n+1}$ which have a constant curvature function $F$ were often called {\it Weingarten hypersurfaces} in the previous literature. It was shown in \cite[Theorem 28]{EspinarGalvezMira:/} that round spheres are the only examples of compact Weingarten hypersurfaces in the hyperbolic space $\mathbb H^{n+1}$. In this spirit, we obtain similar rigidity results for $F$ radially symmetric instead of constant:  
\begin{thm}\label{thm:elliptic}
	Let $n\geq 2$, $K\in \bbR$, $|\alpha| \geq 1$ and $M^{n}\sub \bbM_{K}^{n+1}$ be a convex hypersurface, which is located in the northern hemisphere for $K>0$, such that
	\eq{\label{thm:elliptic-A}\mrm{Sec}_M\geq -\al K.}
	Suppose that $F$ is a convex function satisfying \autoref{cond}, which is a solution to 
	\eq{\label{CP}F = \ga {\rm c}_{_{K}}^{\al}}
	for some constant $\ga$. Then $M$ is a geodesic sphere, which is centred at the origin provided $K\neq 0$. If $\al=1$, the convexity assumption on $M$ can be dropped, while for $\al=-1$ the convexity of $F$ may be replaced by inverse concavity.
\end{thm}

Such results have been obtained for  $\sigma_{\ell}$ by integral methods: for instance, \cite{WuXia:06/2014} studies constant linear combinations of higher order mean curvatures, \cite{Wu:07/2016} analyses constancy of $\co  \sigma_\ell$ in $\mathbb H^{n+1}$, and   \cite{KwongLeePyo:/2018} deals with hypersurfaces having radially symmetric higher order mean curvatures in general $\bbM_{K}^{n+1}$ under mild convexity assumptions. But those integral techniques are restricted to the $\sigma_{\ell}$ because they are divergence free in spaceforms. Our maximum principle approach enables us to relax the assumptions on the curvature functions, at the cost of having to impose a condition on the sectional curvature of the hypersurface. However, notice that if $(1+ \al) K >0$ this assumption is weaker than convexity, while if this product is $\leq 0$ the condition already implies convexity.

\subsection*{Classification of solitons}

In the study of singularity formation along curvature flows, especially the mean curvature flow, the class of {\it{self-shrinking solutions}}, simply called {\it{solitons}} subsequently, plays an important role. For the mean curvature flow in Euclidean space they arise as blow-up limits of type-I singularities (see \cite{Huisken:/1990}), and they satisfy the elliptic equation
\eq{H = \ip{x}{\nu}.}
Huisken \cite{Huisken:/1990} showed that the only compact mean-convex solitons are spheres.

 A similar recent result when $H$ is replaced by the Gauss curvature $K$, see \cite{BrendleChoiDaskalopoulos:/2017}, settled the long standing open problem of whether the flow by certain powers of the Gauss curvature of $n$-dimensional hypersurfaces, $n\geq 3$, converges to a round sphere; the convergence to a soliton had already been proved in \cite{AndrewsGuanNi:/2016}. 
 
 The study of solitons for more general curvature functions has received plenty of attention, also in spaceforms; see e.g. \cite{GaoLiMa:08/2018, GaoLiWang:11/2020, GaoMa:02/2019, McCoy:/2011}. Here one considers the general equation
\eq{ \label{soliton} F^{\be} = u,}
where $\be\in \bbR$, $F$ is a function of the principal curvatures with suitable assumptions and \eq{\label{spt-def} u = {\rm s}_{_{K}}(r) \ip{\del_{r}}{\nu}}
is the generalised support function.
From a well-known duality relation by means of the Gauss map for hypersurfaces of the sphere and itself, respectively a duality relation between hypersurfaces of the hyperbolic and De Sitter space, from \autoref{thm:elliptic} we can deduce a new classification result for convex solitons in the sphere $\mathbb S^{n+1}_1$, as well as in the upper branch of the $(n+1)$-dimensional De Sitter space  $\bbS^{n,1}$ with sectional curvature $K=1$, i.e.
\eq{\bbS^{n,1} = \{y\in \bbR^{n+2}\cn -(y^{0})^{2}+\sum_{i=1}^{n+1}(y^{i})^{2} = 1, y^0 >0 \}.} More precisely, with the notation
\eq{ \label{signat}\sgn(\bbM) = \begin{cases} 1, &\bbM = \bbS^{n+1}_{1}, \\
		-1, &\bbM = \bbS^{n,1},
\end{cases}}
we prove the following result
\begin{cor}\label{cor:elliptic}
	Let $n\geq 2$,  $\abs{\be}\leq 1$, $\be\neq 0$ and $\bbM$ be either $\bbS^{n+1}_{1}$ or $\bbS^{n,1}$.  Consider $M^{n}\sub \bbM$ a closed strictly convex hypersurface, and in case $\frac{1-\beta}{\beta} \sgn(\bbM) >0$, we assume further that
	\eq{\mrm{Sec}_M \leq \frac{\sgn(\bbM)}{1-\beta}.}
If $F$ is an inverse convex function satisfying \autoref{cond}, which is a solution to the soliton equation \eqref{soliton}, then $M$ is a geodesic sphere centred at the origin. In case $\be=1$, the inverse convexity may be replaced by concavity.
\end{cor}

\begin{rem}
	\autoref{cor:elliptic} is remarkable in several ways. Firstly, to our knowledge this is the first such result, where the $\beta$-regime ranges down to zero. This is surprising, as in the Euclidean space, for $F = K^{1/n}$ and $\be\leq n/(n+2)$ the result is false, see \cite{Andrews:09/2000, BrendleChoiDaskalopoulos:/2017}. Note however that $K^{1/n}$ is not inverse convex. Secondly, in all of the previous results of this type, the inverse concavity of $F$ was exploited crucially. The duality approach allows us to deal with a further class of curvature functions, which could not be treated by earlier methods. 
\end{rem}

Notice that, while Weingarten hypersurfaces are known to be geodesic spheres in $\mathbb S^{n,1}$ (see \cite{Roldan:/2022}), we are not aware of rigidity results for solitons in this setting. On the other hand, to have some model examples in mind, $F = |A|$ satisfies the assumptions of \autoref{thm:elliptic} and the harmonic mean curvature is suitable for \autoref{cor:elliptic}.

\subsection*{The problem of extending a non-local flow to curved spaces} As said before,  Huisken \cite{Huisken:/1987} introduced the {\sc vpmcf} of convex hypersurfaces in the Euclidean space, 
\eq{\label{sphere-VPMCF}\partial_t{x}=\br{\mu(t) - H}\nu,}
where the global term is the average mean curvature $\mu =\overline H = \fint_{M_{t}}H$. Taking
$\mu$ as in \eqref{def-mu} for $K = 0$, McCoy \cite{McCoy:02/2004} obtained convergence of convex hypersurfaces in $\mathbb R^{n+1}$ to round spheres under a flow that preserves any quermassintegral (which in the Euclidean case coincide with the mixed volumes, see \autoref{sec:quermass}).

 Huisken already pointed out that an interesting problem is to extend his result to non-Euclidean ambient spaces, with the warning that the generalization will not be straightforward because \eqref{sphere-VPMCF} does not preserve convexity in general Riemannian manifolds, due to terms with an unfavourable sign in the evolution equation of the second fundamental form. In particular, for hypersurfaces in $\bbM_{K}^{n+1}$ the Weingarten matrix $h^{i}_{j}$ evolves according to
 \[(\partial_t - \Delta) h^i_j = (|A|^2 - n K) h^{i}_{j} + 2 K H \delta^i_j - \mu (h^i_\ell h^\ell_j + K \delta^i_j).\]
 
 Notice that for $K <0$ the bad term is $ 2 K H \delta^i_j$, which comes from the background geometry and causes that convexity is not preserved in general. This failure is independent of the non-local nature of the flow; indeed, if we replace convexity by $h$-convexity ($\kappa_1 > |K|$), Miquel and the first author \cite{Cabezas-RivasMiquel:/2007} proved that $h$-convex hypersurfaces can be deformed under \eqref{sphere-VPMCF} to a geodesic sphere; this was extended by Andrews and Wei \cite{AndrewsWei:10/2018} for a class of quermassintegral preserving flows.  The curvature condition was relaxed to positive sectional curvature ($\kappa_1 \kappa_2 > |K|$) by Andrews, Chen and Wei in \cite{AndrewsChenWei:05/2018} in the volume-preserving case.
 
 Notice that the complication for $K >0$ is of completely different nature, since the fatal term is now $-\mu K \delta^i_j$, and thus comes directly from the global term. Indeed, Huisken illustrated this \cite{Huisken:/1987} with an intuitive example: if the flow starts with a convex hypersurface of $\mathbb S^{n+1}$ with a portion $M^*$ $C^2$-close to the equator, then in this region $\ov H \gg H$ and hence $M^\ast$ moves in the outward direction crossing the equator, and thus the evolving hypersurface becomes non-convex.
 
 This obstruction to the preservation of convexity in an ambient sphere supports the claim that the flow \eqref{sphere-VPMCF} is, geometrically, not the most natural generalisation of the same flow in Euclidean case to the spherical ambient space. Indeed,  our alternative flow \eqref{flow-eq} does preserve pinching of the principal curvatures and hence it succeeds in driving any convex initial hypersurface to a geodesic sphere. Notice that  Huisken's example is actually the motivation for the definition of \eqref{flow-eq}, as the effect of multiplying the global term by $\co(r)$ is to slow down the motion as the hypersurfaces approach the equator. 
 
 In short, to extend Huisken's results to the hyperbolic space one needs to strengthen the notion of convexity, whereas for the ambient sphere we propose a different generalization of the flow (notice that \eqref{flow-eq} and \eqref{sphere-VPMCF} coincide for the Euclidean space), which works for convex hypersurfaces.

\subsection*{The isoperimetric nature of the flow}
 In addition, under \eqref{sphere-VPMCF} the surface area is non-increasing, and hence Huisken's theorem provides an alternative proof of the isoperimetric inequality for convex hypersurfaces of $\bbR^{n+1}$. An interesting side-effect of the extra term in \eqref{flow-eq} is that this flow is no longer of isoperimetric nature in the classical sense, because if we choose $\mu$ to preserve enclosed volume, the surface area is no longer decreasing necessarily.

 However, the flow \eqref{flow-eq} with global term chosen to preserve the weighted volume $\int_{\Om_{t}}{\rm c}_{_{K}}$, has decreasing surface area, which suggests that in principle it is the right flow to prove the isoperimetric type inequality
 \eq{\label{isop-w} \int_{\Om_{0}}{\rm c}_{_{K}}\leq \phi(\abs{M_{0}})}
 with equality if and only if $\Om_{0}$ is a ball centred at the origin. Here $\phi$ is a function that gives equality on the slices. This was originally shown in \cite[Proposition~4]{GiraoPinheiro:12/2017} by other means and hence we do not pursue any further investigation in this matter here.

This reinforces the idea that our new flow has a geometric meaning beyond the generalization of Huisken's result, and we hope that in the future some interesting new applications will follow.

\subsection*{Structure of the paper} The contents of this paper are organized as follows. We first introduce in \autoref{sec:notation} the basic notation and evolution equations that ensure that our flow preserves the quermassintegrals, while \autoref{sec:Convexity} gathers new estimates for strictly convex hypersurfaces in the sphere, which may be of independent interest, like a refined outradius bound (\autoref{opt-outC}) or inradius control in terms of pinching (\autoref{bounds-rad}). Then in \autoref{sec:pinching} we prove that the {\it pinching deficit} decreases exponentially under the flow as time evolves, which is the key to get convergence of the evolving hypersurfaces. To achieve upper curvature bounds, we perform a technically intricate process in \autoref{sec:H-bdd}, which includes a delicate iterative changing of origin to ensure an optimal configuration that enables us to gain some uniform bound on the global term for some controlled time interval. This is a novel method, providing an alternative to proving initial value independent curvature bounds after a waiting time. To finish the proof of \autoref{thm:main}, in \autoref{sec:global} we establish long time existence and convergence to a geodesic sphere is done in \autoref{sec:convergence}. Finally, the elliptic results are proved in \autoref{sec:elliptic}.

\section{Notation, Conventions and Preliminary Results}\label{sec:notation}

\subsection*{Hypersurfaces in space forms}

Let $x\cn M\hra \bbM_{K}^{n+1}$ be an embedded smooth hypersurface in a simply connected space form $\bbM_{K}^{n+1}$ enclosing a bounded domain $\Omega$. Then the metric in polar coordinates is given by
\eq{\bar g = dr^{2}+{\rm s}_{_{K}}^{2}(r)\si,}
where $r$ is the radial distance to a fixed point $\mc O \in \bbM_{K}^{n+1}$ and $\si$ is the round metric on $\bbS^{n}$.

The trigonometric functions in \eqref{trig} satisfy the following computational rules:
\[{\rm c}'_{_{K}}= - K \, {\rm s}_{_{K}}, \quad {\rm c}^2_{_{K}} + K\, {\rm s}^2_{_{K}} =1. \]
We will also use the related notation $ {\rm co}_{_{K}}(r) =
{\rm c}_{_K}(r)/{\rm s}_{_K}(r)$.

For the outward pointing unit normal $\nu$, we define the second fundamental form $h=(h_{ij})$ by
\eq{\bar\n_{X}Y = \n_{X}Y - h(X,Y)\nu,}
where $\bar\n$ is the Levi-Civita connection of the metric  $\bar g =\ip{\cdot}{\cdot}$ on $\bbM_{K}^{n+1}$ and $X,Y$ are vector fields on $M$. 
We adopt the summation convention throughout and latin indices indicate components with respect to a coordinate frame $(\del_{i})_{1\leq i\leq n}$ on the domain of the embedding $x$.

If the induced metric on $M$ is denoted by $g$, then we write $\De$ for  its Laplace-Beltrami operator  and define the Weingarten operator $A = (h^{i}_{j})$ via
\eq{h_{ij} = g(A(\del_{i}),\del_{j}) = g_{ik}h^{k}_{j}.}
Recall that the symmetry of $h$ and the {\it Codazzi equations} 
\[\nabla_i h_{jk} = \nabla_j h_{ik}\]
imply that the tensor $\nabla A$ is totally symmetric. Moreover, one can relate the geometry of a hypersurface $M$ with the ambient manifold $\bbM_{K}^{n+1}$ by means of the {\it Gauss equation}
\eq{\label{Gauss} R_{ijk\ell} = h_{ik} h_{j\ell}  - h_{i\ell} h_{j k} + K(g_{ik} g_{j\ell}  - g_{i\ell} g_{j k}).}

On the other hand, if 
\eq{\kappa_1 \leq \cdots \leq \kappa_n}
denote the eigenvalues of the operator $A$, that is, the principal curvatures of $M$, we consider the normalised mean curvatures $H_\ell$ defined as
\[H_\ell = \binom{n}{\ell}^{-1} \sigma_\ell, \qquad \text{with} \qquad \sigma_\ell = \sum_{1 \leq i_1 < \cdots < i_\ell \leq n} \kappa_{i_1} \cdots \kappa_{i_\ell}.\]
In particular, $H_1 = H/n$ and $H_n$ equals the Gauss curvature. We use the convention that $H_0 = 1$. For convex hypersurfaces, these symmetric functions satisfy the {\it Newton-MacLaurin} inequalities \cite{WangXia:/2014}:
\eq{\label{Newton-Mc} H_{\ell -1} H_k \geq H_\ell H_{k-1} \q \text{for} \q 1 \leq k < \ell \leq n.}

We will also use the following {\it Hsiung-Minkowski} identities \cite{GuanLi:/2015}:
\eq{\label{Mink} (\ell + 1) \int_M  u \sigma_{\ell+1} = (n-\ell) \int_M\co \sigma_\ell}
for $\ell = 0, \ldots, n-1.$

Later on, we need control of the support function from below, given that there is a uniform ball enclosed by the evolving domain. Fortunately, for strictly convex domains such control is easy to obtain. We furnish quantities like $r$ and $u$ with a hat, if their reference point is not the origin. The right reference point will then be displayed as a subscript, and in cases where the reference point is clear by context, it is suppressed

\begin{lemma} \label{u-bdd:lem}
Let $\Om\sub \bbM_{K}^{n+1}$ be a strictly convex domain with $p\in \Om$ and $M=\del\Om$. Then the support function 
\eq{\hat u_{p} = {\rm s}_{_{K}}(\hat r_{p})\ip{\del_{\hat r_{p}}}{\nu},}
where $\hat r_{p}$ is the distance to the point $p$, 
\eq{\hat u_{p}\geq \min_{M}\hat u_{p} = \min_{M}{\rm s}_{_{K}}(\hat r_{p})= {\rm s}_{_{K}}(\dist(p,M)).}
\end{lemma}

\pf{
At a global minimum of the support function we have $\n \hat u_{p}=0$. It is well known (cf.~\cite{GuanLi:/2015}) that 
\eq{\label{na-u} \nabla_i \hat u_{p} = h_i^j \nabla_j \Big(\frac{1-\co}{K}\Big) = - h_i^j \frac{\co'}{K} \nabla_j \hat r_p = {\rm s}_{_{K}}(\hat r_{p}) h_i^j \nabla_j \hat r_{p}.}
Accordingly, due to the invertability of $A$, we also have $\n \hat r_{p}= 0$. Hence at such a point, and for $K>0$,
\eq{\hat u_{p} = {\rm s}_{_{K}}(\hat r_{p})\geq \min({\rm s}_{_{K}}(\min_{M}\hat r_{p}),{\rm s}_{_{K}}(\max_{M}\hat r_{p})),}
due to the concavity of ${\rm s}_{_{K}}$ within the interval $[0,\pi/\rt{K}]$. In case $K\leq 0$, ${\rm s}_{_{K}}$ is increasing, so in this case we are done. Now suppose that
\eq{\label{u-bdd:lem-1}{\rm s}_{_{K}}(\max_{M}\hat r_{p})<{\rm s}_{_{K}}(\min_{M}\hat r_{p}).}
Due to the symmetries of the sine function we must then have
\eq{\min_{M}\hat r_{p}>\fr{\pi}{\rt{K}} - \max_{M}\hat r_{p}.}
At a point $\xi\in M$, where $\hat r_{p}$ is maximized, the geodesic which connects $p$ and $\xi$, intersects $M$ in another point, say $\ze\in M$. Due to the convexity of $M$ there holds 
\eq{\fr{\pi}{\rt K}>\dist(\ze,\xi)= \hat r_{p}(\ze) + \hat r_{p}(\xi)\geq \hat r_{p}(\ze) + \fr{\pi}{\rt{K}}-\min_{M}\hat r_{p}\geq \fr{\pi}{\rt{K}},}
a contradiction. Hence \eqref{u-bdd:lem-1} can not be valid and the proof is complete.
}

\subsection*{Mixed volumes and quermassintegrals} \label{sec:quermass}
We define the {\it curvature integrals} or {\it mixed volumes} as follows:
\eq{V_{n-\ell}(\Omega) = \int_M H_\ell \, d V, \qquad \text{for} \quad \ell =0, \ldots, n.}

On the other hand, for any connected domain $\Omega \subset \bbM_{K}^{n+1}$ bounded by a compact hypersurface of class $C^3$, the {\it quermassintegrals} are given by (see \cite{Solanes:/2006} or \cite[chapter 17]{Santalo:book}):
\eq{\label{def-Wk} W_\ell(\Omega)= \frac{(n+1-\ell) \, \omega_{\ell-1} \cdots \omega_0} {(n+1) \, \omega_{n-1} \cdots \omega_{n-\ell}} \int_{\mc L_\ell} \chi(L_\ell \cap \Omega)\, dL_\ell, \quad \ell = 1, \ldots, n.}
Here $\mc L_\ell$ represents the space of $\ell$-dimensional totally geodesic subspaces $L_\ell$ in $\bbM_{K}^{n+1}$, where one can define a natural invariant measure $dL_\ell$, and $\omega_n = |\mathbb S^n|$ is 
the area of $n$-dimensional unit sphere in $\mathbb R^{n+1}$. If $\Omega$ is a convex set, then the function $\chi$ is equal to 1 if $L_\ell \cap \Omega \neq \emptyset$, and 0 otherwise. 

One typically sets 
\eq{W_0(\Omega) =|\Omega| \q \text{and} \q W_{n+1}(\Omega) = \frac{\omega_n}{n+1}.} 
Moreover, using the Cauchy-Crofton formula (cf.~\cite{Santalo:book}) we recover the area of the hypersurface as 
\eq{|\partial \Omega|= (n+1) W_1(\Omega).}
 Accordingly, volume and area-preserving flows can be regarded as particular cases of quermassintegral preserving flows.

Mixed volumes and quermassintegrals are related (see \cite[Proposition 7]{Solanes:/2006}) in a space of constant curvature $\bbM_{K}^{n+1}$ by means of  
\eq{\label{mixed-quer} \frac1{n+1}V_{n-\ell}(\Omega) &= W_{\ell+1}(\Omega) - K \frac{\ell}{n+2-\ell} W_{\ell-1}(\Omega), \qquad \ell =1, \ldots, n, \\
V_n(\Omega) &= (n+1) W_1(\Omega) = |\partial \Omega|.}
Notice that in $\mathbb R^{n+1}$ the mixed volumes coincide with the quermassintegrals, up to a constant factor.

The following result gathers the evolution equations of the above defined quantities under a normal variation.
\begin{lemma} \label{basic-evol}
If $M_t$ is a hypersurface of $\bbM_{K}^{n+1}$ evolving along a flow given by $\partial_t x = \varphi \nu$, then 
\begin{enumerate}
	\item[{\rm (a)}] $\displaystyle \partial_t {\rm Vol}(\Omega_t) = \int_M \varphi \, dV_t\quad$ and 
	 $\quad \displaystyle \partial_t |M_t| = \int_M \varphi \, H \, dV_t$.
	\item[{\rm (b)}] $\displaystyle \partial_t \int_M H_\ell \, dV_t  = \int_M \varphi \big((n-\ell) H_{\ell+1} - K \ell H_{\ell -1}\big) \, dV_t,\q\ell = 0,\dots, n.$
	\smallskip
	\item[{\rm (c)}] $\displaystyle \partial_t W_\ell(\Omega_t)   = \frac{n+1-\ell}{n+1}\int_M \varphi \, H_\ell \, dV_t, \q \ell = 0,\dots, n$.
	\item[{\rm (d)}] $\displaystyle \del_{t}g_{ij} = 2\varphi h_{ij}.$
	\smallskip
	\item[{\rm (e)}] $\displaystyle \del_{t}h^{i}_{j} = -g^{ik} \nabla^2_{kj}\varphi - \varphi h^{i}_{k}h^{k}_{j}-K\varphi \delta^i_j$.
\end{enumerate} 
\end{lemma}

\pf{Formulas in (a) and (b) were deduced in \cite{Reilly:/1973}. The evolution in (c) follows arguing by induction on $\ell$ and using the relation \eqref{mixed-quer}; this was done in \cite[Proposition 3.1]{WangXia:/2014} for $K =-1$. The evolution for the metric and the Weingarten operator are standard, e.g. \cite[chapter 2]{Gerhardt:/2006}.}

\begin{cor}
If the global term in \eqref{flow-eq} is chosen as in \eqref{def-mu}, then the quermassintegral $W_\ell(\Omega_t)$ is constant along the flow \eqref{flow-eq}.
\end{cor}
The fact that $\mu(t) >0$ for strictly convex hypersurfaces is heavily used within the proof of \autoref{thm:main}.

\begin{rem}
Notice that a global term given by 
\[\mu(t) = \frac{\int_M  H \big((n-\ell) H_{\ell+1} - K \ell H_{\ell -1}\big) \, d V_t}{\int_M  {\rm c}_{_{K}} \big((n-\ell) H_{\ell+1} - K \ell H_{\ell -1}\big) \, d V_t}\]
leads to a flow that preserves the mixed volume $V_{n-\ell}(\Omega_t)$. Unlike the quermassintegral-preserving case, this term does not have a sign for convex hypersurfaces if $K >0$. For $K < 0$ this difficulty disappears, but another type of mixed volume preserving curvature flows for $h$-convex hypersurfaces in the hyperbolic space was already studied in \cite{Makowski:08/2012}.
\end{rem}

We use the following conventions for the use of constants. Indexed letters $C$, i.e.~ $C_{0}, C_{1}$, etc. will retain a specific meaning throughout the whole paper, while the letter $C$ denotes a generic constant, which is always allowed to change from line to line and depends on the quantities listed in the formulation of the lemma or theorem. Capital letters also stand for ``large'' constants. A similar convention holds for lower case letters, which stand for ``small'' constants. The only exception from this convention concerns the use of various versions of the letter $t$, like $T$, $\tau$, $\hat\tau$ etc. Those always refer somehow to time and $t$ denotes the time variable, while $T$, $\tau$, $\hat\tau$ etc. will, once defined, not change value.

\section{Geometry and location of strictly convex hypersurfaces in the sphere}\label{sec:Convexity}

This section presents some geometric results for strictly convex hypersurfaces of the sphere, which are required to prove \autoref{thm:main}. In particular, we obtain inradius estimates in terms of pinching, as well as a suitable outball configuration in terms of pinching and the value of any given $W_{\ell}(\Om)$. {\bf{Throughout \autoref{sec:Convexity}, \autoref{sec:pinching} and \autoref{sec:H-bdd}}, we make the standing assumption that $M\sub \bbS^{n+1}_{K}$ is a strictly convex hypersurface enclosing a domain $\Om$.}

Let $B_r$ denote a geodesic ball of radius $r$ in $\bbS_{K}^{n+1}$. The outer radius of $\Om$ is given by 
\eq{\rho_+(\Om) = \inf\big\{R >0 \, |\, \Omega \subset B_R(q) \text{ for some } q \in \bbS_{K}^{n+1}\big\},}
 and the inner radius is given by
\eq{\rho_-(\Om) = \sup\big\{\rho >0 \, |\,  B_\rho(p) \subset \Omega \text{ for some } p \in \bbS_{K}^{n+1}\big\}.}
 In \cite{Andrews:/1994b}, it was shown that pinched hypersurfaces of the Euclidean space satisfy a uniform control of outer by inner  radius, and a version for a positive ambient space can be found in \cite[Sec.~6]{Gerhardt:/2015}. The proof of this version relied on uniform positivity of the smallest principle curvature, which is insufficient for our purposes. Hence we provide a more general version in the following proposition.
 
\begin{prop} \label{in-out}
 If for some number $C_{0}>0$ there holds the pinching estimate $\kappa_n \leq C_0 \kappa_1$ in $M$, then the outer radius is estimated from above according to
	\eq{\rho_+(\Om) \leq C_1 \, \rho_-(\Om),}
for some positive constant $C_1=C_{1}(n,K,C_{0})$.
\end{prop}

\pf{For simplicity but without loss of generality, we assume $K = 1$. Due to a classical result \cite{CarmoWarner:/1970},
	\eq{0<\rho_{+}(\Om)< \fr{\pi}{2},}
	because $M$ lies in some open hemisphere. Hence there is a centre $q$ such that
	\eq{\Om\sub B_{\rho_{+}}(q),}
	and it is true that $q\in \bar\Om$ (cf.~\cite[p.~455]{Santalo}). By moving $q$ slightly inwards, we can achieve $M\sub B_{\fr{\pi}{2}}(q)$ and it is starshaped around $q$.
	
	Now consider the stereographic projection from the antipodal point $-q$, where $q$ is mapped to the origin $0 \in \bbR^{n+1}$. It follows  (see \cite[(6.15)]{Gerhardt:/2015}) that the metric $\bar g$ of $\mathbb S^{n+1}$ is conformal to the Euclidean metric; more precisely,
	\eq{\bar g =  e^{2\psi}(dr^{2}+r^{2}\si), \q \text{with} \q \psi(r) = - \ln(1 + r^2/4).}
 Hereafter, we denote by tilde the Euclidean geometric quantities. On $
 \bar B_{\fr{\pi}{2}}(q)$ the metric $\bar g$ is uniformly equivalent to the Euclidean metric.

 Next, from \cite[(1.1.51)]{Gerhardt:/2006} we get
	\eq{e^{\psi}\ka_{i} = \ti\ka_{i} + d \psi(\ti\nu)}
	and hence from our pinching assumption,
	\eq{0<C_{0}^{-1}\leq \fr{\ka_{1}}{\ka_{n}} = \fr{\ti\ka_{1} + d \psi(\ti\nu)}{\ti\ka_{n}+ d \psi(\ti\nu)}\leq 1.}
 Thus
	\eq{\ti\ka_{1} + d \psi(\ti\nu)\geq C_{0}^{-1}(\ti\ka_{n} + d \psi(\ti\nu))}
	and
	\eq{\ti\ka_{1}\geq C_{0}^{-1}\ti\ka_{n} +(C_{0}^{-1}-1) \psi'(r)\ip{\partial_r}{\ti\nu}\geq C_{0}^{-1}\,\ti\ka_{n},}
	because $\psi'<0$ and $\ti M$ is starshaped, i.e.~$\ip{\del_{r}}{\ti\nu}>0$. Therefore the Euclidean hypersurface $\ti M\sub \ti B_{2}(0)$ is pinched, which from \cite[Lemma~5.4]{Andrews:/1994b} leads to
	\eq{\rho_{+}(\Om)  \leq C\ti\rho_{+}(\Om) \leq C \ti\rho_{-}(\Om)\leq  C_1 \rho_{-}(\Om),}
	where we have used the uniform equivalence of the ambient metrics.
}

\begin{cor} \label{bounds-rad}
If for some number $C_{0}>0$ there holds the pinching estimate $\kappa_n \leq C_0 \kappa_1$ in $M$, then one can find positive constants $d_{1}$ and $C_{2}$, depending on $n$, $K$, $C_{0}$ and $W_\ell(\Omega)$, such that
	\eq{ d_{1} \leq \rho_{-}(\Om)\leq C_{2}<\fr{\pi}{2\rt{K}}.}
\end{cor}

\pf{By the definition of inner and outer radius, we can find points $p, q \in \mathbb \bbS^{n+1}_{K}$ such that
	\eq{B_{\rho_{-}(\Om)}(p) \subset \Omega \subset B_{\rho_{+}(\Om)}(q).}

	From \eqref{def-Wk}, the quermassintegrals $W_\ell$ are clearly monotone under the inclusion of convex domains, and hence
	\eq{W_{\ell}\big(B_{\rho_{-}(\Om)}(p)\big) \leq W_{\ell}(\Omega) \leq W_{\ell}\big(B_{\rho_{+}(\Om)}(q)\big).}
We obtain with \autoref{in-out},
	\eq{C_{2}:=f_{\ell}^{-1}(W_{\ell}(\Om))\geq \rho_{-}(\Om)\geq C_{1}^{-1}\rho_+(\Om) \geq C_{1}^{-1} f_{\ell}^{-1}(W_\ell(\Omega))=:d_{1},}
 where $f_{\ell}$ denotes the increasing function given by  $f_{\ell} (r) = W_{\ell}(B_r)$. 
 }
 
 The trivial outer radius estimate
 \eq{\rho_{+}(\Om)<\fr{\pi}{2\rt{K}}}
 is not good enough for our purposes. Now we present a refined estimate which should be of independent interest in the future.

\begin{thm} \label{opt-outC}
	The outer radius satisfies
	\eq{\rho_{+}(\Om)\leq \fr{\pi}{2\sqrt{K}}-\fr{\log W_{\ell}(\cH)-\log W_{\ell}(\Om)}{(n+1-\ell)\max_{M}H}=:\fr{\pi}{2\rt{K}} - \fr{d_{2}}{\max_{M}H}, }
	where $\cH$ is an open hemisphere.
\end{thm}

\pf{
	From the initial hypersurface $M_{0}=M$, we start the curvature flow
	\eq{\label{ICF} x\cn [0,T^{*})\x \bbS^{n} &\ra \bbS^{n+1}_K\\
		\partial_t x &= \fr{H_{\ell-1}}{H_{\ell}}\nu.}
	Then \autoref{basic-evol} (c) ensures that $W_{\ell}$ evolves in time according to
	\eq{\del_{t}W_{\ell}(\Om_{t}) = \fr{n+1-\ell}{n+1} \, V_{n-\ell+1}(\Om_{t})\leq (n+1-\ell)W_{\ell}(\Om_{t}),}
	where the inequality follows from \eqref{mixed-quer}.
	
	From \cite{Gerhardt:/2015,MakowskiScheuer:11/2016}, we know that \eqref{ICF} preserves the strict convexity and the solution converges smoothly to an equator, while we also have the estimate
	\eq{W_{\ell}(\Om_{t})\leq W_{\ell}(\Om)e^{(n+1-\ell)t}.}
	Hence the maximal existence time of \eqref{ICF} is at least
	\eq{T^{*}\geq \fr{1}{n+1-\ell}\br{\log W_{\ell}(\cH)-\log W_{\ell}(\Om)},}
	because we know that $W_{\ell}(\Om)$ must converge to $W_{\ell}(\cH)$ at $T^{*}$. 
	
	Now \cite[Lemma~4.7]{MakowskiScheuer:11/2016} leads to the curvature bound
	\eq{\max_{M_{t}}H\leq \max_{M_{0}}H \q \text{for all}\q 0\leq t<T^{*}. }
	Due to the convexity this implies a full second fundamental form bound, as well as a bound
	\eq{\label{quot-bound} \fr{H_{\ell}}{H_{\ell-1}}\leq H_1 \leq \max_{M_{0}}H,}
	which follows by application of \eqref{Newton-Mc} for $k = 1$.
	
	Next let $E=\del\cH$ be the limiting equator of the flow and let $r$ be the radial distance from the centre of $\cH$, which contains all $M_{t}$. 
	Define
	\eq{\ti r(t) = \max_{\bbS^{n}}r(t,\cdot) = r(t,\xi_{t}),}
	where $\xi_{t}$ is any choice of point where the maximum is realised. 
	The function $\ti r$ is Lipschitz and hence differentiable a.e. At times, where $\ti r$ is differentiable, we have
	\eq{\fr{d}{dt}\ti r(t) = \fr{H_{\ell-1}}{H_{\ell}}(t,\xi_{t}),}
	where we used that $\nu(t,\xi_{t})=\del_{r}$.
	Integration and \eqref{quot-bound} yield
	\eq{\fr{\pi}{2 \sqrt{K}} - \ti r(0) = \ti r(T^{*}) - \ti r(0) = \int_{0}^{T^{*}}\fr{d}{dt}\ti r\geq \fr{\log W_{\ell}(\cH)-\log W_{\ell}(\Om)}{(n+1-\ell)\max_{M_{0}}H}.}
	Hence
	\eq{\max_{M_{0}}r=\ti r(0)\leq \fr{\pi}{2 \sqrt{K}} - \fr{\log W_{\ell}(\cH)-\log W_{\ell}(\Om)}{(n+1-\ell)\max_{M_{0}}H}.}
	Accordingly, $M_{0}$ fits into a neighbourhood of the origin of size given by the right hand side of the latter inequality, and therefore the outer radius is controlled by the very same quantity.
}

\begin{rem}
A lune, i.e. the intersection of two hemispheres, shows that an outer radius bound in terms of $W_{\ell}(\Om)$ can not be independent of $\max H$.
\end{rem}

\autoref{opt-outC} enables us to find, for a given strictly convex hypersurface of the sphere, a suitable origin, which allows a ball of controlled size within $\Om$ and at the same time ensures a controlled positive distance of $M$ to the equator.

\begin{lemma}\label{config}
	 There exists an origin $\cO\in \Om$, such that with the constant $d_{2}$ from \autoref{opt-outC} there holds
	\eq{\label{config-A}B_{4\ep}(\cO)\sub \Om\q\mbox{and}\q \max_{M}r\leq \fr{\pi}{2 \sqrt{K}} - 4\ep}
	for all $\displaystyle {\ep\leq \fr 14\min\br{\fr{d_{2}}{2\max_{M}H},\fr{\fr{\pi}{2} - \tan^{-1}(\max_{M}H/\rt{K})}{2\rt{K}}}.}$
\end{lemma}

\pf{
If $B_{\rho_{+}(\Om)}$ denotes an outball for $\Om$ with center $\cO$, then we know that $\cO\in \bar\Om$. The distance of focal points from $M$ can be calculated from the evolution of the Weingarten operator in \autoref{basic-evol} along the normal variation with speed $\vp=-1$. Then the largest principal curvature is controlled by the solution to the ODE
\eq{y' &= y^{2}+K\\
	y(0)&=\max_{M}H,}
which exists for all $t < t_0$, with
\eq{t_{0} := \fr{\pi}{2\rt{K}} - \fr{\tan^{-1}(\max_{M}H/\rt{K})}{\rt{K}}.}
Hence, around all points belonging to the set $\{x\in \Om\cn \dist(x,M) = t_{0}/2\}$, there exists an interior ball of radius $t_{0}/2$. In addition, if we shift $\cO$ by a distance of $d_{2}/(2\max_{M}H)$ in any direction, there still is the same amount of space between $M$ and the new equator. Therefore, if we shift $\cO$ into $\Om$ along a perpendicular geodesic only by the amount
\eq{\ep\leq\ep_{0} = \fr 14\min\br{\fr{d_{2}}{2\max_{M}H},\fr{t_{0}}{2}},}
then \eqref{config-A} holds.
}

\section{Monotonicity of the Pinching Deficit}\label{sec:pinching}

The geometric results from \autoref{sec:Convexity} depend on the quality of the pinching and the size of the quermassintegral. In the following we investigate how these quantities behave under the flow \eqref{flow-eq}. As this flow is defined to be quermassintegral preserving, the key ingredient for proving \autoref{thm:main} is the pinching estimate to be proven in this section.
 
Again, we assume $K>0$, unless stated otherwise. As a prior step,  we need the following evolution equations.

\begin{lemma} \label{evol}
For every choice of origin $\cO$, along \eqref{flow-eq} the induced metric $g$ and second fundamental form $h$ satisfy the following evolution equations: 
\eq{\partial_t g_{ij} &= 2(\mu {\rm c}_{_{K}}-H) h_{ij}\\
\partial_t h_{j}^{i} &=  \De h_{j}^{i}+ (\abs{A}^{2} - n K) h_{j}^{i} +2KH\de^{i}_{j} - \mu \big( {\rm c}_{_{K}}h_{k}^{i}h^{k}_{j} + K u \, h_{j}^{i}\big),
}
where $u$ is the generalised support function in \eqref{spt-def} of $M_{t}$ with respect to the origin $\cO$.
\end{lemma}

\pf{
The evolution of the metric comes from \autoref{basic-evol} (d).
For the evolution of $A$, we depart from the standard evolution equation
\eq{\label{evolh1}\partial_t h_{ij} = \nabla^2_{ij}(H-\mu {\rm c}_{_{K}}) + (\mu {\rm c}_{_{K}} - H)(h_{ik}h^{k}_{j} - Kg_{ij}),}
see \cite[Theorem 3-15]{Andrews:/1994c} and use the following Simon's type identity:
\eq{\label{lem:ev-h-2}\n^{2}_{ij} H&= \Delta h_{ij}+ \big(\abs{A}^{2}-nK\big)h_{ij}+ H\big(K g_{ij} -h_{i}^{k}h_{jk}\big).}
Now we expand the second derivatives of ${\rm c}_{_{K}}$:
\eq{\label{Hess-c}-\n^{2}_{ij}{\rm c}_{_{K}} =  d{\rm c}_{_{K}}(\nu)h_{ij}-\bar\n^{2}{\rm c}_{_{K}}(x_{i},x_{j}) =  {\rm c}_{_{K}}'\ip{\del_{r}}{\nu}h_{ij} + K{\rm c}_{_{K}} g_{ij},}
where we have used 
\eq{\bar\n^{2}{\rm c}_{_{K}} = -K{\rm c}_{_{K}} \bar g.}
The proof is complete, using the evolution of the metric to revert to $h^{i}_{j}$.
}

\begin{rem} \label{pres-conv}
Notice that the strong maximum principle for tensors applied to the evolution of $h^i_j$ already implies that the property of strict convexity is preserved for all times. Accordingly, $H >0$ and the quotient $\frac{\kappa_1}{H}$ is well-defined as long as the flow exists.
\end{rem}

Next we deduce an evolution equation that is the key to convergence of the flow.

\begin{prop}\label{ev-p-parabolic}
Let $\ka_{1}$ be the smallest eigenvalue of $A$. Then for every choice of origin $\cO$, under the flow \eqref{flow-eq} with initial data $M$, the function
\eq{p = \fr{\ka_{1}}{H}}
is a supersolution to the following evolution equation in viscosity sense, as long as the flow exists:
\eq{\del_{t}p - \De p &= \fr{2}{H}\sum_{k=1}^{n} \sum_{j>D}\fr{(\n_{k}h^{1}_{j})^{2}}{\ka_{j}-\ka_{1}}+2d p(\n\log  H)+\fr{\mu}{H} {\rm c}_{_{K}}\abs{\mr{A}}^{2}p\\
		&\hp{=}+\mu {\rm c}_{_{K}}\ka_{1}(\tfr 1n  -p) + 2nK(\tfr 1n - p),
}
where $D$ is the multiplicity of $\ka_{1}$.
\end{prop}

\pf{
Assume that the flow is defined on a maximal time interval $[0, T)$. Let $(t_{0}, \xi_{0})\in   (0,T) \x M$ and let $\eta$ be a smooth lower support of $p$ at $(t_{0}, \xi_{0})$, i.e. $\eta$ is defined on a spacetime neighbourhood $\cU$ of $(t_{0}, \xi_{0})$ and there holds
\eq{\eta(t_{0}, \xi_{0}) = p(t_{0}, \xi_{0}),\q \eta\leq p_{|\cU}. }
Hence $\vp = H\eta $ is a smooth lower support for $\ka_{1}$. 

Now we take coordinates with the properties
\eq{g_{ij} = \de_{ij},\q h^{i}_{j} = \ka_{j}\de^{i}_{j} \q \text{at} \q (t_{0}, \xi_{0}).}
If we denote by $D$ the multiplicity of $\ka_{1}(t_{0}, \xi_{0})$, then
 at the point $(t_{0}, \xi_{0})$ and for all $1\leq i,j\leq D$ it is satisfied (see \cite[Lemma~5]{BrendleChoiDaskalopoulos:/2017}) that
\eq{\del_{t}h^{i}_{j} = \de^{i}_{j}\del_{t}\vp,\q \n_{k}h^{i}_{j} = \de^{i}_{j}\n_{k}\vp}
and
\eq{\label{Lem5} \n_{kk}^{2}\vp \leq \n^{2}_{kk}h^{1}_{1}-2\sum_{j>D}\fr{(\n_{k}h^{1}_{j})^{2}}{\ka_{j}-\ka_{1}}.}

Next, write 
\eq{G^{i}_{j} = \del_{t}h^{i}_{j} - \Delta h^{i}_{j},}
and compute at $(t_{0}, \xi_{0})$
\eq{\del_{t}\eta &=\fr{\del_{t}h^{1}_{1}}{H}-\fr{\vp}{H^{2}}\del_{t}H\\
			&=\fr{G^{1}_{1}+ g^{kl}\n^{2}_{kl}h^{1}_{1}}{H}-\fr{\vp}{H^{2}}\del_{t}H\\
			&\geq \fr{1}{H}\sum_{k=1}^{n}\br{\n^{2}_{kk}\vp + 2\sum_{j>D}\fr{(\n_{k}h^{1}_{j})^{2}}{\ka_{j}-\ka_{1}}}+\fr{G^{1}_{1}}{H}-\fr{\vp}{H^{2}}\del_{t}H\\
			&= \fr{2}{H}\sum_{k=1}^{n} \sum_{j>D}\fr{(\n_{k}h^{1}_{j})^{2}}{\ka_{j}-\ka_{1}}+\fr{G^{1}_{1}}{H}+\De \eta +2d\eta(\n\log H)-\fr{\eta}{H}G^{k}_{k}.}
			
On the other hand, we know by \autoref{evol} that
\eq{G^{i}_{j} &= \big(\abs{A}^{2}-nK\big)h_{j}^{i} +2KH\de^{i}_{j}-\mu \big({\rm c}_{_{K}} h_{k}^{i}h^{k}_{j} + K u h_{j}^{i}\big).} 
Accordingly, we get
\eq{G^{1}_{1} - \eta G^{k}_{k}&=  -\mu {\rm c}_{_{K}}\ka_{1}^{2} +\eta\mu {\rm c}_{_{K}}\abs{A}^{2}+ 2KH(1-n\eta) . }
Finally, by means of
\eq{\abs{A}^{2} = \abs{\mr{A}}^{2}+\tfr{1}{n}H^{2},}
we obtain
\eq{G^{1}_{1} - \eta G^{k}_{k}&=  \eta\mu {\rm c}_{_{K}}\abs{\mr{A}}^{2}+\mu {\rm c}_{_{K}}\ka_{1}H(\tfr 1n  -\eta)+ 2KH(1-n\eta), }
which completes the proof.
}

\begin{cor}\label{pinching}
For every choice of origin $\cO$, for which $M\sub \cH(\cO)$, the flow \eqref{flow-eq} with initial data $M$ stays in $\cH(\cO)$ and it improves every pinching, i.e. the {\it{pinching deficit}}
\eq{ \om(t) =\fr{1}{n} - \min_{M_{t}}\fr{\ka_{1}}{H}}
is exponentially decreasing,
\eq{\om(t)\leq \om(s)e^{-2nK(t-s)}\q\text{for all} \q 0\leq s\leq t<T,}
where $T$ is the maximal time of existence of the flow with initial data $M$.
\end{cor}

\pf{
We need the evolution equation of ${\rm c}_{_{K}}$. From \eqref{flow-eq} and \eqref{Hess-c} we obtain
\begin{align}
 \del_{t}{\rm c}_{_{K}} &= {\rm c}_{_{K}}' \partial_t r = -K {\rm s}_{_{K}} dr(\partial_t x)  = K u(H - \mu {\rm c}_{_{K}}) \nonumber = \De {\rm c}_{_{K}} + K{\rm c}_{_{K}} (n - \mu u). \label{evol_cK}
\end{align}
The preservation of ${\rm c}_{_{K}}>0$ follows immediately by strong maximum principle, as long as the flow exists.

The strict convexity is also preserved from \autoref{pres-conv}. The statement about the pinching deficit follows from the strong maximum principle for viscosity solutions, e.g. see \cite{Da-Lio:09/2004}, and from the fact that
\eq{\del_{t}\om \leq -2nK\om,}
where we can discard the terms including $\mu$ because ${\rm c}_{_{K}}>0$ and $\mu(t) >0$ by convexity.
}

In particular, for $s = 0$ we reach a pinching relation between the biggest and smallest principal curvatures of $M_t$:
\eq{\label{pinching-k}\kappa_1 \geq \Big(\min_{M_0} \frac{\kappa_1}{H}\Big) H \geq C_0^{-1} \kappa_n,}
provided that an origin is chosen such that the strictly convex initial hypersurface is contained in the open hemisphere centred at that origin.

\section{Upper Curvature Bounds}\label{sec:H-bdd}

Notice that, unlike in previous treatments of quermassintegral preserving curvature flows, an upper bound for the global term does not come automatically from an upper bound for $H$, since the $\co$ in the denominator of $\mu$ is not uniformly bounded away from zero, at least not without further work. Moreover, we need some uniform control of $\mu$ to get bounds for $H$.
 
 To overcome these difficulties, the idea is to choose the origin such that a configuration as in \autoref{config} is achieved, which will allow us to deduce uniform bounds on the curvature and the global term in a short but controlled interval $[0, \tau(\ep)]$ (\autoref{Hbdd-lem}). Then, since the pinching is at least as good as at the beginning, we can repeat this process as often as needed, in order to keep the flow going as long as we like (\autoref{lem:global}). During this evolution, the pinching improves exponentially and at some point will be so strong, that $M_{t}$ is very close to a sphere. From here the flow is very easy to estimate and no further shifting of the origin is necessary. Now we implement all required steps to make this argument rigorous.

 A key idea to obtain a curvature bound is to adapt a well-known trick from \cite{Tso:/1985}, which consists in a suitable combination of the generalised support function with the mean curvature. For this we need control on the size of inballs during the flow.

\subsection*{A lower bound for the support function of an arbitrary inball}
For a domain $\Om$ and a point $p\in \Om$, we say that $B$ is an {\it{inball at $p$}}, if $p$ is the center of $B$ and $B$ has maximal radius with the property that $B\sub \Om$.
In the sequel we are going to prove, that along the flow, the radii of inballs at $p$ don't decrease too quickly. We give a quantitative estimate.

We need to be careful because now we are dealing with two different support functions: we denote by $u$ the support function with respect to the origin  $\cO$ that is implicit in the flow equation, while for a given point $p\in \Om$, 
\eq{\hat u\equiv\hat u_{p} = {\rm s}_{_{K}} (\hat r_{ p}) \ip{\partial_{\hat r_{p}}}{\nu}} takes another interior point $p$ as origin of distances. Accordingly, $r$ and $\hat r$ mean distance from $\cO$ and $p$, respectively. Similar notations will apply to the corresponding trigonometric functions, i.e.
\eq{\hat {\rm s}_{_{K}}:={\rm s}_{_{K}}(\hat r_{p})\q\mbox{and}\q \hat{\rm c}_{_{K}}:={\rm c}_{_{K}}(\hat r_{p}).} Note that for brevity we suppressed the dependence on the point $p$ within the notations $\hat u$ and $\hat r$.

\begin{lemma} \label{evolHu}
For every choice of origin, along \eqref{flow-eq} the evolution equations of the mean curvature $H$ and the support function $\hat u=\hat u_{p}$  are given by
\eq{\partial_t H & = \Delta H + H (|A|^2 + K n) - \mu (\co |A|^2 + u K H), \smallskip \\
	\partial_t \hat u & = \Delta \hat u + \hat u |A|^2 + (\mu \co  - 2 H) \hat{\rm c}_{_K} + \mu  K \s \hat{\rm s}_{_K} \big<\nabla r, \nabla \hat r\big>.}
\end{lemma}

\pf{The formula for $H$ follows directly by taking the trace in the evolution equation for $h^{i}_{j}$ from  \autoref{evol}.
On the other hand, a standard calculation leads to
\eq{\partial_t \hat u & = \ip{\bar \nabla_t(\hat{\rm s}_{_K} \partial_{\hat r})}{\nu} + \ip{\hat{\rm s}_{_K} \partial_{\hat r}}{ \bar \nabla_t\nu} \\
& = (\mu \co - H)\hat{\rm c}_{_K} + \hat{\rm s}_{_K}\ip{\partial_{\hat r}}{\nabla(H-\mu \co)} \\
&= \Delta \hat u + (\mu \co - 2 H) \hat{\rm c}_{_K} + \hat u |A|^2 - \mu \,\hat{\rm s}_{_K} \ip{\partial_{\hat r}}{\nabla \co}, }
where we applied well-known formulas for $\Delta \hat u$ and $\bar \nabla_t(\hat{\rm s}_{_K} \partial_{\hat r})$ on $\bbM_{K}^{n+1}$ (cf.~\cite[(4.6) and (4.11)]{Cabezas-RivasMiquel:/2007}). The stated formula follows by realizing that \eq{\ip{\partial_{\hat r}}{ \nabla \co} = - K \s \ip{\partial_{\hat r}}{\nabla r}.}
}

\begin{prop} \label{comp-pple}
For every choice of origin $\cO$, for which $M\sub \cH(\cO)$, and for every $p\in \Om$ and radius $\rho$ with the property $B_{\rho}(p)\sub\Om$, the solution $M_t = \partial \Omega_t$ of \eqref{flow-eq}, with initial data $M$ and maximal existence time $T>0$, satisfies the following:
\enum{
\item There is a positive constant $\ti\tau=\ti\tau(n,K,\rho)$ with the property
\eq{B_{ \rho/4}(p) \subset \Omega_{t} \quad \text{ for all } \q
t \in [0,\min(\ti\tau,T)).}
\item One can find positive constants $d_{3}\leq \frac1{2\sqrt{K}}$ and $\tau$, depending on $n,K,C_{0}$ and $W_{\ell}(\Om)$, with the property
	\eq{\hat u_{p_{\Om}} - 2d_{3} \geq  2d_{3} > 0 \quad \text{ for all } \q
t \in [0,\min(\tau,T)),}
where $p_{\Om}$ is the center of an inball corresponding to the inradius $\rho_{-}(\Om)$.
}

\end{prop}
\pf{
(i)~Let us first obtain the evolution of the distance $\hat r =\hat r_{p}$ from the fixed point $p$ to the points on $M_t$ under the flow \eqref{flow-eq}:
\eq{\label{pf:comp-pple-1}\partial_t \hat r = d \hat r (\partial_t x)   = (\mu {\rm c}_{_{K}}-H) \ip{\nu}{\partial_{\hat r}}.}

On the other hand, $\rr(t)$ denotes the radius of a geodesic sphere centred at $p$ that moves under the ordinary mean curvature flow starting at $\rr(0)= \rho/2$, that is, 
\eq{ \rr'(t) = - n \, {\rm co}_{_{K}}(\rr(t)),}
whose solution is given by 
\eq{{\rm c}_{_{K}}(\rr(t)) = e^{K nt} {\rm c}_{_{K}}(\rho/2) \q \text{for} \q t \geq 0. }
As $\co$ is a decreasing function, 
\eq{\rr(t) \geq  \rho/4 \q \text{ if and only if } \q e^{K nt} {\rm c}_{_{K}}(\rho/2) \leq {\rm c}_{_{K}}(\rho/4),}
meaning that
\eq{\rr(t) \geq \rho/4 \quad \text{if and only if} \quad t \leq \frac{1}{K n} \log \frac{{\rm c}_{_{K}}(\rho/4)}{{\rm c}_{_{K}}(\rho/2)}=:\ti\tau.}

 Set $f(t,\cdot) = \hat r(t,\cdot) -\rr(t)$ for $t \in [0,\min(\tau,T))$. Then $f(0,\cdot)>0$ and $f$ evolves according to 
\eq{\partial_t f = (\mu {\rm c}_{_{K}}-H) \ip{\nu}{\partial_{\hat r}} + n \, {\rm co}_{_{K}}(\rr(t)). }
If there exists a first time $t_1$ such that the geodesic sphere $B_{\rr(t_1)}$ touches the hypersurface $M_{t_1}$ at some point $x_1$, then at this first minimum for $f$ it holds $H(x_1, t_1) \leq n {\rm co}_{_{K}}(\rr(t_1))$, $\ip{\partial_{\hat r}}{\nu}=1$ and $\partial_t f(x_1, t_1) \leq 0$. Consequently, taking into account that $\co >0$ and strict convexity is preserved, we have
\[\partial_t f (x_1, t_1) \geq \mu(t_1) \, {\rm c}_{_{K}} (r(t_1)) >0, \]
which is a contradiction, and hence the statement follows.

(ii) Apply (i) with $p=p_{\Om}$ and $\rho=\rho_{-}(\Om)$ and obtain the desired $\tau$ due to \autoref{bounds-rad}. Using \autoref{u-bdd:lem} we get
\eq{\label{d3} \hat u_{p_{\Om}} \geq {\rm s}_{_{K}} (\rho_{-}(\Om)/4)  \geq {\rm s}_{_{K}} (d_{1}/4) =: 4 d_{3},}
where we applied the lower bound in \autoref{bounds-rad}. 
}

\subsection*{An upper bound for the mean curvature}

To estimate the mean curvature along a solution with suitably located initial data, we use the well known auxiliary function 
\eq{\Phi_{p} = \fr{H}{\hat u_{p}-d_{3}},}
 which, after choosing the origin $\cO$ as in \autoref{comp-pple}, is well-defined for a while for some suitable $p\in \Om$. Routine computations lead to the evolution of $\Phi$ where we suppress the dependence on $p$ within the notation.

\begin{lemma} \label{evol-aux}

Under the assumptions of \autoref{comp-pple}, along the flow \eqref{flow-eq}
the function $\Phi$ evolves according to
\eq{\partial_t \Phi  &= \Delta \Phi + \frac{2}{\hat u - d_{3}} \ip{\nabla \Phi}{\nabla  \hat u} 
+ \Phi \br{n K -\frac{d_{3}}{\hat u -d_{3}} |A|^2} + 2 \Phi^2   \hat{\rm c}_{_K}  \\
	  &\hp{=}  - \mu K \Phi u -  \frac{\mu}{\hat u-d_{3}} \Big(\co |A|^2 + \Phi(K \s  \hat{\rm s}_{_K} \ip{\nabla r}{\nabla \hat r} + \co  \hat{\rm c}_{_K})\Big).}
\end{lemma}

Unlike in previous literature, we cannot neglect all the terms including $\mu$, as some of them do not have a sign. Nor is it known at this point that $\mu$ is bounded. The novelty about our approach is to make use of \autoref{config}, a configuration which enables us to gain some control on $\mu$ and then get an estimate on $\Phi$ and $H$. 

Another complication arises from the necessity of using an iterative change of origin. The configuration of \autoref{config} depends on curvature. Hence we need a very precise estimate of curvature as the flow progresses and it is insufficient to estimate the curvature by a multiple of its initial value, as then our time interval, along which $B_{\ep}(\cO)\sub \Om_{t}$ is valid, would decrease and we would not be able to prove long-time existence. 

For this reason, we introduce a novel method, which also gives an interesting alternative to proving initial value independent curvature bounds after a waiting time as for example in \cite[Equation~(17)]{McCoy:02/2004}. It provides a bound on $\Phi$, which is uniform in $p$ lying within a certain region. 

With this purpose, we define a modified auxiliary function 
\eq{\Psi\cn[0,T)\x \bbS^{n}\x \bbS^{n+1}_{K}&\ra \bbR\\
						(t,\xi,p)&\mt \begin{cases}(\dist(p,M_{t})-2d_{3})\Phi_{p}, & p\in V_{t}\\
							0	,&~\mbox{else}, 
						\end{cases}\\
					}
with
\eq{V_{t}:=\{p\in \Om_{t}\cn \dist(p,M_{t})>2d_{3}\},}
and where for $p\notin \bar\Om$ the distance to $M_{t}$ is defined to be negative. In short,
\eq{\Psi(t, \xi, p)= \max(0,\min_{M_{t}}\hat r_{p}-2d_{3})\Phi_{p}.}

Note that $\Psi$ is Lipschitz, because when $\min_{M_{t}}\hat r_{p} = \dist(p,M_{t})\geq 2d_{3}$, then by \autoref{u-bdd:lem}
\eq{\hat u_{p}\geq {\rm s}_{_{K}}(\min_{M_{t}}\hat r_{p})\geq {\rm s}_{_{K}}(2d_{3}),}
due to $\min_{M_{t}}\hat r_{p}<\pi/(2\rt{K}).$ Furthermore,
\eq{\hat u_{p}\geq {\rm s}_{_{K}}(2d_{3})=\tfr{1}{\rt K}\sin(2d_{3}\rt{K})\geq \tfr{6}{5}d_{3},}
where we used $\sin x\geq \tfr 35 x$ for $x\in [0,\pi/2]$. Hence $\Phi_p$ is well-defined for $p \in \bar V_t$. In the sequel we write, for brevity,
\eq{\max_{M_{t}}\Psi = \max_{\bbS^{n}\x\bbS^{n+1}_{K}}\Psi(t,\cdot).}

\begin{lemma} \label{Hbdd-lem}
There exists an origin $\cO\in \Om$, a constant $C_{3}(n,K,C_{0},W_{\ell}(\Om))$, and constants $\tau$ and $\ep_{1}$, depending on $n,K,C_{0},W_{\ell}(\Om)$, $\max(C_{3},\max_{M_{0}}\Psi)$, such that for the solution $M_{t}=\del \Om_{t}$ of \eqref{flow-eq} with initial data $M_{0}=M$ and maximal existence time $T>0$ there hold:
\enum{
	\item $M_{t}\sub \cH(\cO)$ for all $t\in [0,\min(\tau,T)]$.
	\item $ \hat u_{p_{\Om}}\geq 4d_{3}$ and  $B_{\ep_{1}}(\cO)\sub\Om_{t}$ for all $t\in [0,\min(\tau,T)).$ \medskip
	
	\item	The function $\Psi$ satisfies
	\eq{\max(C_{3},\max_{M_{t}}\Psi)\leq \max(C_{3},\max_{M}\Psi)\q\mbox{for all}\q t\in [0,\min(\tau,T)).}	
}
\end{lemma}

\pf{
Define the functions
\eq{\ep_{0}(y) = \fr 14\min \br{\fr{d_{2}d_{3}\rt{K}}{2y},\fr{\fr{\pi}{2}-\tan^{-1}(\fr{y}{d_{3}K})}{2\rt{K}}},\q \text{for }y>0,}
and 
\eq{q(y) = \fr{ n\pi C_0^{\ell +1}}{\rt{K}}\br{1+\fr{5\pi}{\rt{K}d_{3}}}\fr{y}{\ep_{0}(y)}+\fr{n\pi^{2}}{4}y+\fr{\pi}{\rt{K}}y^{2}-\fr{d_{3}^{2}}{n}y^{3},}
where $d_{3}$ and $C_{0}$ are the constants from \autoref{comp-pple} and \eqref{pinching-k} respectively. It is clear that $q(y)$ converges to $-\8$ as $y\ra \8$ and hence it has a largest zero $\bar y$, which only depends on $n,K,C_{0}$ and $W_{\ell}(\Om)$. Let $\psi_{0}=\max_{t=0}\Psi$ and define
\eq{ \label{def-eps1} \ep_{1} = \ep_{0}(\max(\bar y,\psi_{0})).}

(i)~\&~(ii):~With $p_{\Om}$ from \autoref{comp-pple} (ii) and by definition of $d_{3}$ in \eqref{d3}, along $M_{0}=M$ there holds
\eq{H = (\hat u_{p_{\Om}}-d_{3})\Phi_{p_{\Om}}\leq \fr{\Phi_{p_{\Om}}}{\rt{K}}\leq \fr{{\rm s}_{_{K}}(\min_{M}\hat r_{p_{\Om}})-2d_{3}}{2d_{3}}\fr{\Phi_{p_{\Om}}}{\rt{K}}\leq \fr{1}{d_{3}\rt{K}}\psi_{0}}
and hence
\eq{\ep_{1}\leq\ep_{0}(d_{3}\rt{K}\max_{M}H)<\fr{\pi}{2\rt{K}},}
where the latter estimate is due to the definition of $\ep_{0}$. This exactly the threshold required to apply \autoref{config}.

Then we can apply \autoref{config} with $\ep=\ep_{1}$ in order to obtain a suitable origin $\cO\in\Om$ with the property \eqref{config-A}. 
From the first part of \autoref{comp-pple} applied to $p=\cO$ and $\rho = 4\ep_{1}$, as well as from the second part of \autoref{comp-pple}, we obtain $\tau = \tau(n,K,C_{0},W_{\ell}(\Om),\ep_{1})$, up to which the claimed properties of (ii) are satisfied. Property (i) is then clear from the fact that at the equator ${\rm c}_{_{K}} = 0$ and hence $\tfr{d}{dt}\max r<0$.

(iii)~ 
Next, we bound the function $\Psi$. Suppose $\Psi$ attains a positive maximum over the set $[0,\bar t]\x\bbS^{n}\x \bbS^{n+1}_{K}$ at some $(\bar t,\bar\xi,\bar p)$. Define the Lipschitz function
\eq{\bar\Psi_{\bar p}(t)= \max_{\bbS^{n}}\Psi(t,\cdot,\bar p),}
which is positive in some small interval $J:=[\bar t-\de,\bar t]$.
Thus we have 
\eq{\dist(\bar p,M_{t})>2d_{3}\q\fa t\in J.}
Hence in $J$ the function $\Phi_{\bar p}$ is smooth and 
\eq{\bar\Psi_{\bar p}(t) = (\min_{M_{t}}\hat r_{\bar p}-2d_{3})\max_{M_{t}}\Phi_{\bar p}}
is differentiable almost everywhere in $J$.

There holds for almost every $t\in J$,  using \eqref{pf:comp-pple-1} and $\ip{\del_{\hat r_{\bar p}}}{\nu}>0$ due to $\bar p\in \Om_{ t}$ for all $t\in J$,
\eq{\fr{d}{dt}\bar\Psi_{\bar p}&= \fr{d}{dt}\min_{M_{t}}\hat r_{\bar p}\max_{M_{t}}\Phi_{\bar p} + (\min_{M_{t}}\hat r_{\bar p}-2d_{3})\fr{d}{dt}\max_{M_{t}}\Phi_{\bar p}\\
						&=({\rm c}_{_{K}}(r)\mu - H(\bar t,x(\bar t,\mrm{argmin}_{M_{\bar t}}\hat r_{\bar p})))\ip{\del_{\hat r_{\bar p}}}{\nu}\max_{M_{t}}\Phi_{\bar p}\\
						&\hp{=}+ (\min_{M_{\bar t}}\hat r_{\bar p}-2d_{3})\fr{d}{dt}\max_{M_{t}}\Phi_{\bar p}\\
						&\leq \mu\max_{M_{t}}\Phi_{\bar p}+ (\min_{M_{\bar t}}\hat r_{\bar p}-2d_{3})\fr{d}{dt}\max_{M_{t}}\Phi_{\bar p}.
}
Inserting the evolution equation from \autoref{evol-aux}, discarding two good terms, and using $\n\Phi_{\bar p}( t,\xi_{t})=0$ at all maximisers $\xi_{t}$, we obtain
\eq{\fr{d}{dt}\bar\Psi_{\bar p}	&\leq \mu\max_{M_{t}}\Phi_{\bar p}+ (\min_{M_{\bar t}}\hat r_{\bar p}-2d_{3})\left( \max_{M_{t}}\Phi_{\bar p} \br{n K -\tfrac{d_{3}}{\hat u_{\bar p} -d_{3}} |A|^2} + 2 \max_{M_{t}}\Phi_{\bar p}^2  \right.  \\
	  &\hp{=} \left. -  \frac{\mu}{\hat u_{\bar p}-d_{3}} \Big(\max_{M_{t}}\Phi_{\bar p}(K \s  \hat{\rm s}_{_K} \ip{\nabla r}{\nabla \hat r_{\bar p}} + \co  \hat{\rm c}_{_K})\Big)\right),
}
where we used that $u$ and ${\rm c}_{_{K}}$ are positive due to $M_{ t}\sub \cH(\cO)$ and $\cO\in \Om_{ t}$. 
 Using ${\rm s}_{_{K}}\leq K^{-1/2}$ and
\eq{\fr{\min_{M_{\bar t}}\hat r_{\bar p}-2d_{3}}{\hat u_{\bar p}-d_{3}}(K{\rm s}_{_{K}}\hat{\rm s}_{_{K}}\ip{\n r}{\n \hat r_{\bar p}}+{\rm c}_{_{K}}{\hat{\rm c}_{_{K}}})\leq 2\cdot\fr{\pi/(2\rt K)}{\tfr 15 d_{3}} = \fr{5\pi}{d_{3}\rt{K}},}
we get
\eq{\fr{d}{dt}\bar\Psi_{\bar p}	 &\leq \br{1+\fr{5\pi}{\rt{K}d_{3}}}\mu\max_{M_{t}}\Phi_{\bar p}+nK\bar\Psi_{\bar p}+2\bar\Psi_{\bar p}\max_{M_{t}}\Phi_{\bar p}-\fr{d_{3}^{2}}{n}\bar\Psi_{\bar p}\max_{M_{t}}\Phi_{\bar p}^{2}.
}
By means of the Hsiung-Minkowski identity \eqref{Mink} and \autoref{u-bdd:lem}, we estimate for all $0\leq t\leq \tau$,
\eq{\label{bo-mu} \mu(t) \leq n\binom{n}{\ell+1} \fr{\int_{M_{t}}\ka_{n}^{\ell+1}}{\int_{M_{t}}u\si_{\ell+1}}\leq \frac{n}{{\rm s}_{_{K}}(\ep_{1})} \frac{\int_{M_t} \kappa_n^{\ell +1}}{\int_{M_t} \kappa_1^{\ell +1}} \leq \frac{2 n C_0^{\ell +1}}{\ep_{1}},  }
where the constant $C_0$ comes from the pinching \eqref{pinching-k}.
Hence
\eq{\fr{d}{dt}\bar\Psi_{\bar p}  &\leq 2 n C_0^{\ell +1}\br{1+\fr{5\pi}{\rt{K}d_{3}}}\fr{\max_{M_{t}}\Phi_{\bar p}}{\ep_{0}(\max(\bar y,\psi_{0}))}+nK\bar\Psi_{\bar p}\\
		&\hp{=}+2\bar\Psi_{\bar p}\max_{M_{t}}\Phi_{\bar p}-\fr{d_{3}^{2}}{n}\bar\Psi_{\bar p}\max_{M_{t}}\Phi_{\bar p}^{2}\\
	&= \fr{2 n C_0^{\ell +1}}{\min_{M_{\bar t}}\hat r_{\bar p}-2d_{3}}\br{1+\fr{5\pi}{\rt{K}d_{3}}}\fr{\bar\Psi_{\bar p}}{\ep_{0}(\max(\bar y,\psi_{0}))}+nK\bar\Psi_{\bar p}\\
	&\hp{=}+\fr{2}{\min_{M_{\bar t}}\hat r_{\bar p}-2d_{3}}\bar\Psi_{\bar p}^{2}-\fr{d_{3}^{2}}{n(\min_{M_{\bar t}}\hat r_{\bar p}-2d_{3})^{2}}\bar\Psi_{\bar p}^{3}.
}

Multiplication with $(\min_{M_{\bar t}}\hat r_{\bar p}-2d_{3})^{2}$ gives for almost every $t\in J$,
\eq{(\min_{M_{\bar t}}\hat r_{\bar p}-2d_{3})^{2}\fr{d}{dt}\bar\Psi_{\bar p}	&\leq \fr{ n\pi C_0^{\ell +1}}{\rt{K}}\br{1+\fr{5\pi}{\rt{K}d_{3}}}\fr{\bar\Psi_{\bar p}}{\ep_{0}(\max(\bar y,\psi_{0}))}\\
				&\hp{=}+\fr{n\pi^{2}}{4}\bar\Psi_{\bar p}+\fr{\pi}{\rt{K}}\bar\Psi_{\bar p}^{2}-\fr{d_{3}^{2}}{n}\bar\Psi_{\bar p}^{3}.
}
which is strictly negative whenever $\max(\bar y,\psi_0)< \bar\Psi_{\bar p}(t)$. To conclude the argument, suppose that 
\eq{\Psi(\bar t,\bar\xi,\bar p)=\max_{[0,\bar t]\x\bbS^{n}\x\bbS^{n+1}_{K}} \Psi>\max(\bar y,\psi_{0}).}
We know that 
\eq{\bar\Psi_{\bar p}(t)\leq \Psi(\bar t,\bar\xi,\bar p)}
with equality at $t=\bar t$. But before we have shown that, for $t$ close to $\bar t$, we have $\tfr{d}{dt}\bar\Psi_{\bar p}<0$ almost everywhere, which is impossible.
Hence we obtain the desired estimate with $C_{3} = \bar y$.
}

\subsection*{Higher order curvature bounds}

We use the estimate from \autoref{Hbdd-lem} to control the global term and estimate the derivatives of curvature.

\begin{lemma}\label{A-der-bdd}
For the origin $\cO\in \Om$ from \autoref{Hbdd-lem} and the solution $M_{t}=\del\Om_{t}$ of \eqref{flow-eq} with initial data $M$ and maximal existence time $T>0$ and for all $m\in \bbN$, there exists $C_{4} = C_{4}(n,m,K,C_{0},W_{\ell}(\Om),\max(C_{3},\max_{M}\Psi))$ with the property
\eq{\mu(t)+\abs{\n^{m}A}\leq C_{4}\q\fa t\in [0,\min(\tau,T)).}
where $\tau$ is the number from \autoref{Hbdd-lem}.
In particular there holds $T\geq\tau$ and the flow exists smoothly on $[0,\tau]$. 
\end{lemma}

\pf{
Up to the time $\min(\tau,T)$, \eqref{bo-mu} holds, i.e. with the notation from the proof of \autoref{Hbdd-lem} we have
\eq{\mu(t)\leq \fr{2nC_{0}^{\ell+1}}{\ep_{1}} = \fr{2nC_{0}^{\ell+1}}{\ep_{0}(\max(C_{3},\psi_{0}))}.}
The curvature derivative bound can be proved by a well-known induction argument, as for example in \cite{Huisken:/1984}. First, due to convexity and \autoref{Hbdd-lem} (ii) and (iii),
\eq{\abs{A}^{2}\leq H^{2}\leq \fr{\Phi^{2}_{p_{\Om}}}{K}\leq \fr{1}{d_{3}^{2}K}(\max(C_{3},\max_{M}\Psi))^2.}
Assuming that all derivatives up to order $m-1$ are bounded by a constant of the form $C_{4}$, we obtain the evolution equation of $\abs{\n^{m}A}$,

\eq{\partial_t |\nabla^m A|^2 & \leq \Delta |\nabla^m
 A|^2- 2 |\nabla^{m+1}
 A|^2 + C (\mu u + K) |\nabla^m
 A|^2
\\ & +  C\!\!\! \sum_{i+j+k = m} \Big(\mu |\nabla^i \co| + |\nabla^i A|\Big) |\nabla^j A| |\nabla^k A||\nabla^m A|,}
where we used that $\nabla u = A \ast \nabla \co$ (see \eqref{na-u}). Here $S \ast T$ denotes any linear combination of tensors formed by contracting $S$ and $T$ by means of $g$. 

Then we claim 
\eq{ \label{DA-bdd}
|\nabla^m A|^2 \le C\q\fa t\in [0,\min(\tau,T)).}
 As $\mu$ and $u$ are bounded, we get
\eq{(\partial_t - \Delta)|\nabla^m A|^2 & \leq - 2 |\nabla^{m+1}
	A|^2 + C|\nabla^m
	A|^2
 +  C\Big( \sum_{i= 0}^m |\nabla^i \co| +1\Big)|\nabla^m A|.}

Notice that $\co(r)$ and $|\nabla  \co|$ are bounded as well. Moreover, from \eqref{Hess-c}, one has for $\ell \geq 0$ the following covariant derivatives
\eq{\nabla^{\ell+2} \co = \nabla^\ell \co * K + u \ast \nabla^\ell A + \sum_{i+j+k = \ell} \nabla^i \co \ast \nabla^j A \ast \nabla^k A,}
which are controlled by uniform constants arguing by induction. In short, we reach 
\eq{(\partial_t - \Delta)|\nabla^m A|^2 & \leq - 2 |\nabla^{m+1} A|^2 +  C|\nabla^m A|\br{1
	+  |\nabla^m A|},}
which leads to \eqref{DA-bdd} by standard maximum principle arguments, as for example in the proof of \cite[Theorem 4.1]{Huisken:/1987}.

As the right hand side of the flow equation and all higher derivatives of the curvature remain uniformly bounded, we conclude (as in \cite[pages 257 ff.] {Huisken:/1984})
that, if $T < \tau$, then $M_t$ converges (as $t \to T$, in the $C^\infty$-topology) to a unique, smooth and strictly convex hypersurface\footnote{Note that due to the bound on $\mu$, it can be seen from \autoref{evolHu} that $H$ is bounded from below on every finite time interval.}. Now we
can apply  short time existence to continue the solution after $T$,
contradicting the maximality of $T$. Hence the solution of \eqref{flow-eq} starting at a strictly convex hypersurface exists on $[0, \tau)$. On this interval we have uniform smooth estimates, and hence the flow also exists on $[0,\tau]$.
}

\section{Construction of a global solution} \label{sec:global}
In the previous section we achieved existence and uniform estimates of any solution to \eqref{flow-eq} with strictly convex initial data $M$ on a time interval $[0,\tau]$ from \autoref{Hbdd-lem}, the length of which only depends on preserved data of the problem. Those are in particular the hemisphere $\cH(\cO)$, the pinching constant $C_{0}$, the quermassintegral $W_{\ell}(\Om)$ and the number $\max(C_{3},\max_{M}\Psi)$. The full curvature derivative bounds also only depend on those quantities. 

Hence we can start an iteration process and shift, at time $i\tau$ with $i\in \bbN$, the origin according to \autoref{Hbdd-lem} applied to the new strictly convex initial hypersurface $M_{i\tau}$. The constant $C_{4}$ from \autoref{A-der-bdd} is then uniform among the integers $i$, because it only depends on quantities which are always preserved. The following lemma makes this precise.

\begin{lemma}\label{lem:global}
Let $M_{0}\sub \bbS^{n+1}_{K}$ be a strictly convex hypersurface enclosing a domain $\Om$. Then there exists a sequence of origins $(\cO_{i})_{i\in \bbN\cup\{0\}}$ and positive numbers $\tau_{0},\ep_{1}(0)$ depending only on $n,K,C_{0},W_{\ell}(\Om_{0}),\max(C_{3},\max_{M_{0}}\Psi)$, such that the problem
\eq{\label{lem:global-A}\del_{t}x &= (\mu_{i}(t){\rm c}_{_{K}}(r_{i})-H)\nu,\q \fa t\in [i\tau_{0},(i+1)\tau_{0})\\
		x(0,\bbS^{n}) &= M_{0}\\
		x((i+1)\tau_{0},\bbS^{n})& = \lim_{t\nearrow (i+1)\tau_{0}}M_{t}, }
where $r_{i}$ is the distance to $\cO_{i}$ and $\mu_{i}$ is given as in \eqref{def-mu} to keep the quermassintegral $W_\ell(\Omega_t)$ for any $\ell=0, 1, \ldots, n$ fixed, has a solution 
\eq{x\cn [0,\8)\x \bbS^{n}\ra \bbS_{K}^{n+1}.} For every $t\geq 0$, $M_{t}$ is strictly convex and there holds
\eq{\label{lem:global-B}B_{\ep_{1}(0)}(\cO_{i})\sub \Om_{t}\q\mbox{and}\q M_{t}\sub \cH(\cO_{i})\q\fa t\in[i\tau_{0},(i+1)\tau_{0}).}The mappings $x(t,\cdot)$ satisfy spatial $C^{\8}$-estimates which are uniform in time.
\end{lemma}

\pf{

For $M_{0}$, pick $\cO_{0}$ according to \autoref{Hbdd-lem}. From \autoref{A-der-bdd} we conclude that the solution $M_{t}$ of \eqref{flow-eq} with initial data $M_{0}$ exists on $[0,\tau_{0}]$, where $\tau_{0}$ and $\ep_{1}(0)$ depend on $n,K,C_{0},W_{\ell}(\Om_{0}),\max(C_{3},\max_{M_{0}}\Psi)$.
The derivatives of $A$ are bounded by 
\eq{C_{4}=C_{4}\br{n,m,K,C_{0},W_{\ell}(\Om_{0}),\max(C_{3},\max_{M_{0}}\Psi)}.}

Now suppose that for $i\geq 0$ the hypersurface $M_{i\tau_{0}}$, the origin $\cO_{i}$ and the solution $(M_{t})_{t\in[i\tau_{0},(i+1)\tau_{0})}$ were already constructed, such that \eqref{lem:global-B},
\eq{\max(C_{3},\max_{M_{t}}\Psi)\leq \max(C_{3},\max_{M_{0}}\Psi),} as well as 
\eq{\label{pf:global-1}\mu_{i}(t)+\abs{\n^{m} A}\leq C_{4}\br{n,m,K,C_{0},W_{\ell}(\Om_{0}),\max(C_{3},\max_{M_{0}}\Psi)}}
all hold for all $t\in [i\tau_{0},(i+1)\tau_{0}]$.
Then apply \autoref{Hbdd-lem} to the initial hypersurface $M_{(i+1)\tau_{0}}$ and obtain an origin $\cO_{i+1}$, such that the solution $M_{t}$ of \eqref{flow-eq} with initial data $M_{(i+1)\tau_{0}}$ satisfies \eqref{lem:global-B},
\eq{\max(C_{3},\max_{M_{t}}\Psi)\leq \max(C_{3},\max_{M_{(i+1)\tau_{0}}}\Psi)\leq \max(C_{3},\max_{M_{0}}\Psi)}
and 
\eq{\mu_{i+1}(t)+\abs{\n^{m} A}\leq C_{4}\br{n,m,K,C_{0},W_{\ell}(\Om_{0}),\max(C_{3},\max_{M_{0}}\Psi)}}
during the interval $[(i+1)\tau_{0},(i+2)\tau_{0}]$ and with the same $\ep_{1}(0)$. Here we also used that $C_{0}$ and $W_{\ell}$ are preserved. This means that the construction and be carried out infinitely often to obtain the desired long-time solution.  
}

\section{Asymptotic estimates and convergence to a spherical cap} \label{sec:convergence}

In the previous sections we have put ourselves into a position where we have a strictly convex flow $(M_{t})_{0\leq t<\8}$ in the sphere. This flow is not necessarily smooth in time, but it satisfies spatial $C^{k}$-estimates which are uniform with respect to time and it has a uniformly bounded global term, due to the proof of \autoref{lem:global}.

Additionally, by means of \autoref{pinching} and the curvature bounds, we get 
\eq{
 \sum_{i=1}^{n}(\ka_{i}-\ka_{1})= H-n\ka_{1}&
		\leq nH\om(t)\leq n C e^{-2nKt},}
which implies exponential decay of the traceless second fundamental form:
\eq{\label{traceless} |\mr{A}|^2 = \frac1{n} \sum_{i <j} (\kappa_j - \kappa_i)^2 \leq C e^{-4nKt}.}

Using this property, in the following we are going to apply some recent estimates of almost-umbilical type due to De Rosa and Gioffr\'e \cite{De-RosaGioffre:/2021}, to show that the process of picking new origins actually terminates after finitely many steps and that the flow will then converge to a geodesic sphere of a uniquely determined radius. The crucial ingredient is the following result.

\begin{thm}\cite[Thm.~1.3]{De-RosaGioffre:/2021}\label{DGR}
Let $n\geq 2$, let $\Si$ be a closed hypersurface in $\mathbb R^{n+1}$ and let $p>n$ be
given. We assume that there exists $c_{0}>0$ such that $\Si$ satisfies the conditions 
\eq{|\Si| = |\bbS^{n}|,\q \|A\|_{L^{p}(\Si)}\leq c_{0}.}
There exist positive numbers $\de,C>0$, depending only on $n,p,c_{0}$, with the following property: if
\eq{\|\mr{A}\|_{L^{p}(\Si)}\leq \de,}
then there exists a vector $c = c(\Si)$, such that $\Si - c$ is a graph over the sphere, namely there
exists a parametrization
\eq{\psi\cn \bbS^{n} \ra \Si,\q \psi(x) = e^{f(x)}x+c,}
 and $f$ satisfies the estimate
 \eq{\|f\|_{W^{2,p}(\bbS^{n})}\leq C\|\mr{A}\|_{L^{p}(\Si)}.}
\end{thm}

In the following we will use this result to prove that the surfaces become exponentially $C^{2}$-close to geodesic spheres and that the necessity to pick new origins vanishes. 

\begin{lemma}\label{DGR-Cor-Conf}
In the situation of \autoref{lem:global}, there exists an integer $m>0$ depending on $n,K$ and $M$, such that in \autoref{lem:global} the origins $\cO_{i}$, $i>m$, may be chosen constantly equal to $\cO_{m}$. 
\end{lemma}

\pf{
Let $m$ be a positive integer to be specified during the proof. 
Let $\mc O_m$ be the flow origin associated to the interval $I_{m}:=[m \tau_{0},(m+1)\tau_{0})$. By stereographic projection from the antipodal point of $\mc O_m$, the family $(M_{t})_{t\in I_{m}}$ can be viewed as a flow in the Euclidean space, then denoted by $(\ti M_{t})_{t\in I_{m}}$. Geometric quantities of this flow, denoted by a tilde, are related to the original ones as follows, see \cite[Prop.~1.1.11]{Gerhardt:/2006}, where $e^{2\vp}$ is the conformal factor:
\eq{g = e^{2\vp}\ti g,\q \nu = e^{-\vp}\ti\nu}
\eq{e^{\vp} A = \ti A + d\vp(\ti\nu)\id, \q e^{\vp} H = \tilde H + n \, d\vp(\ti\nu). }
In particular, we obtain
\eq{\mr{\ti A} = e^{\vp}\mr{A}.}
The surface areas of $M_{t}$ and $\ti M_{t}$ are related by
\eq{\abs{M_{t}} = \int_{M_{t}} 1 = \int_{\ti M_{t}} e^{n\vp}}
and hence
\eq{C^{-1}\abs{M_{t}}\leq \abs{\ti M_{t}}\leq \abs{M_{t}},}
where $C$ depends on $\abs{\vp}_{C^{0}(M_{t})}$.
Now define the scaled hypersurface
\eq{\hat M_{t} = \la\ti M_{t},\q \text{with} \q \la^n = \fr{\abs{\bbS^{n}}}{\abs{\ti M_{t}}},}
so that $\abs{\hat M_{t}} = \abs{\bbS^{n}}$. Now the associated Weingarten operator is
\eq{\hat A = \la^{-1}\ti A = \fr{e^{\vp}}{\la}A - \la^{-1}d\vp(\ti\nu)\id}
and similarly for the traceless Weingarten operator. Hence
\eq{\|\hat A\|_{L^{\8}(\hat M_{t})}\leq C(\|A\|_{L^{\8}(M_{t})}+1)}
and
\eq{\|\mr{\hat A}\|_{L^{\8}(\hat M_{t})}\leq C\|\mr{A}\|_{L^{\8}(M_{t})},}
where $C$ depends on $n$, $\abs{\vp}_{C^{1}(M_{t})}$ and $\abs{M_{t}}$.

Then $\|\hat A\|_{L^{\8}(\hat M_{t})}$ is bounded by \autoref{A-der-bdd} and \eqref{traceless} ensures that $\|\mr{\hat A}\|_{L^{\8}}$ is as small as needed for $m$ big enough. Therefore, we can apply \autoref{DGR} for sufficiently large $m$, to get a function $\hat f$ which, from the embedding theorems of Sobolev spaces
into H\"older spaces, satisfies
\eq{\|\hat f\|_{C^{1}(\bbS^{n})}\leq C\|\mr{\hat A}\|_{L^{\8}(\hat M_{t})}\leq C \|\mr A\|_{L^{\8}(M_{t})}\leq Ce^{-2nK t}\q \text{for all}\q t\in I_{m}.}

Then, due to our curvature bounds \eqref{DA-bdd}, we have full $C^{k}$-bounds on $\hat f$ for all $k$. By iteration of interpolation arguments for $C^k$ bounds (see \cite[Corollary 6.2]{Gerhardt:/2011}), this implies that
\eq{\|\hat f\|_{C^{k}(\bbS^{n})}\leq Ce^{-2nK t}\q\text{for all}\q t\in I_{m}.}
In other words, $\hat M_{t}$ is exponentially $C^{k}$-close to a sphere $\hat\cS_{t}$ for all $k\in \bbN$ and for all $t\in I_{m}$. As the area along the $M_{t}$ is uniformly bounded above and below by \autoref{bounds-rad}, we get $C^k$ bounds for the conformal factor as well, and this property of closeness to a sphere $\cS_{t}$ carries over to the original flow in $\bbS^{n+1}_{K}$. Note that the radii of the spheres $\cS_{t}$ converge to a well defined limit, which ist strictly less than $\pi/(2\rt{K})$, determined by the initial value of $W_{\ell}(\Om_{0})$. Hence the curvature of $M_{t}$ is uniformly bounded from below.

On the other hand, the radial distance to the origin $\mc O_m$ satisfies 
\eq{\label{radial_eq}\partial_t{r} = (\mu \co - H)v^{-1}, \qquad \text{with} \qquad
	v^{2} = 1+\s^{-2}\abs{dr}_{\si}^{2}.} 
Hence, for an error $\de_{m}$ that converges to zero when $m\ra \8$,
\eq{\label{pf:main-0}\partial_t{r} = \br{\fr{\int_{M_{t}}\si_{\ell}H}{\int_{M_{t}}\si_{\ell}{\rm c}_{_{K}}}{\rm c}_{_{K}}-H}v^{-1} \leq \br{\fr{{\rm c}_{_{K}}}{\fint_{\cS_{t}}{\rm c}_{_{K}}}-1}\fr{H_{\cS_{t}}}{v} + \de_{m},}
where $H_{\cS_{t}}$ is the mean curvature of the sphere $\cS_{t}$.
At points which maximize $r$, ${\rm c}_{_{K}}$ is minimized. At such points the first term on the RHS of \eqref{pf:main-0} is strictly negative if ${\rm c}_{_{K}}$ is not constant. Hence, for large $m$, if $\cS_{t}$ is uniformly off-center, $\max_{M_{t}}r$ is decreasing and a similar estimate shows that $\min_{M_{t}}r$ is increasing. Hence, from then on, there is no need to adjust the origin anymore. 
}

\subsection*{Convergence to a spherical cap}

To complete the proof of \Autoref{thm:main}, it only remains to show that the immortal solution coming from \autoref{lem:global} actually converges to a limit geodesic sphere. After $m$-fold picking of a new origin, we now may without loss of generality assume that origins have not been changed at all. We will exploit the $C^{\8}$-estimates for the flow hypersurfaces $M_{t}$ coming from \eqref{DA-bdd}. 
We already know from \autoref{DGR-Cor-Conf}, that every limit point of the flow must be round sphere.

Now we prove that only the sphere centred at the origin can arise as a limit. Notice that the radius of $R$ of any limit sphere is determined by the initial hypersurfaces by means of the equality $W_\ell(B_R) = W_\ell(\Omega_0)$. Denote by $H_{R}$ the mean curvature of such a sphere $S_{R}$. 
Hence, from the evolution \eqref{radial_eq} of the radial distance and for an error $\de$ that converges to zero when $M_{t_{k}}\ra S_{R}$, we get
\eq{\label{pf:main-1}\partial_t{r} =  \br{\fr{{\rm c}_{_{K}}}{\fint_{M}{\rm c}_{_{K}}}-1}\fr{H_{R}}{v} + \de.}
As above, at points which maximize $r$, ${\rm c}_{_{K}}$ is minimized, thus at such points the first term on the RHS of \eqref{pf:main-1} is strictly negative if ${\rm c}_{_{K}}$ is not constant. Therefore the function $\max_{M_{t}}r$ is strictly decreasing in sufficiently small $C^{2}$-neighbourhoods of any non-centred sphere, which excludes those as limits. Thus subsequential limits are unique and the whole flow must converge. 

\section{The elliptic case: rigidity results} \label{sec:elliptic}

In order to prove \autoref{thm:elliptic}, we first have to get the elliptic viscosity equation of the pinching deficit for general curvature function $F$. Let us first gather some prerequisites about these functions.

\subsection*{Symmetric curvature functions} \label{symmetric}
 If $M$ is a hypersurface of $\bbM_{K}^{n+1}$, then we set
\eq{F(x) = f(\ka_{1}(x),\dots,\ka_{n}(x)),}
which can be alternatively seen as a function defined on the diagonalisable endomorphisms, $F=F(A)$, or as a function of a symmetric and a positive definite bilinear form, $F = F(g,h)$. In the latter case, we write
\eq{F^{ij} = \fr{\del F}{\del h_{ij}},\q F^{ij,kl} = \fr{\del^{2}F}{\del h_{ij}\del h_{kl}}.} 
With these conventions, the covariant derivatives are given by
\eq{\label{derF} \nabla_i F = F^{jk} \nabla_i h_{jk}.}

On strictly convex hypersurfaces $M$, we can define the so-called {\it{inverse curvature function}} by
\eq{\ti F(A) = \fr{1}{F(A^{-1})}.}
A curvature function $F$ is called {\it{inverse concave}}, if $\ti F$ is concave. Notice that concavity/convexity with respect to the matrix variables is equivalent to the same property with respect to the eigenvalues (see \cite{Andrews:/2007, Gerhardt:/2006} for deeper information). 
 Next we gather several useful properties:
\begin{lemma} \label{invF}

(a) If $F$ is inverse concave and $M$ is strictly convex, then 
\eq{\label{inv-conc} F^{ij,kl}\n_{1}h_{ij}\n_{1}h_{kl}+ 2\sum_{j}\fr{F^{ii}}{\ka_{j}}(\n_{1}h_{ij})^{2}\geq \fr{2}{F}(\nabla_1 F)^{2} .}

(b)  Under \autoref{cond}, it holds:
\begin{itemize}
\item[(i)] If $F$ is convex, then it is inverse concave (cf.~\cite[Lemma~2.2.12, Lemma~2.2.14]{Gerhardt:/2006})	
\item[(ii)] Euler's formula $F^{ii} \kappa_i = F$, implies that $F$ is strictly positive.
\end{itemize}
\end{lemma}

\subsection*{Rigidity for radial curvature functions}

We start with a result that contains an elliptic version of \autoref{ev-p-parabolic}.
\begin{lemma} \label{lem-el}
Suppose $F$ satisfies \autoref{cond}, then:
\enum{
\item
The Weingarten operator satisfies the following elliptic equation:
\eq{\label{ev:elliptic-h}- F^{rs}\n^{2}_{rs}h_{ij} &=F^{pq,rs}\n_{i}h_{pq}\n_{j}h_{rs} - \n_{ij}^2F \\
		&\hp{=}+F^{rs}h_{ms}h^{m}_{r}h_{ij}+KFg_{ij}-Fh_{mj}h^{m}_{i}-KF^{rs}g_{rs}h_{ij}. }
\item For the function $p = \fr{\ka_{1}}{F},$ it holds
\eq{- F F^{kl}\n^{2}_{kl}p &\geq 2 \sum_{j>D}\fr{F^{kk}(\n_{k}h^{1}_{j})^{2}}{\ka_{j}-\ka_{1}}+F^{pq,rs}\n_{1}h_{pq}\n_{1}h_{rs} +2F^{kl}\n_{k}p\n_{l} F\\
			&\hp{=} + p F^{kl}\n_{kl}F -\n_{11}F
		 +F^{ii}{\rm Sec}_{i1}(\ka_{i}-\ka_{1})}
in viscosity sense, where $D$ is the multiplicity of $\kappa_1$.
}
\end{lemma}

\pf{

(i) We differentiate \eqref{derF}, to get
\[\n_{ij}^2F = F^{pq,rs}\n_{i} h_{pq}\n_{j}h_{rs} + F^{rs}\n^{2}_{ij}h_{rs}.\]
Now combine the commutator formula for second covariant derivatives with Gauss equation \eqref{Gauss} to deduce:
\eq{F^{rs}\n^{2}_{ij}h_{rs} &=  F^{rs}(\n^{2}_{rs}h_{ij} + R^{m}_{isj}h_{mr} + R^{m}_{rsj}h_{mi})\\
		&=F^{rs}\n^{2}_{rs}h_{ij}+F^{rs}(h_{ms}h_{ij}-h_{mj}h_{is}+Kg_{ms}g_{ij}-Kg_{mj}g_{is})h_{r}^{m} \\
		&\hp{=}+ F^{rs}(h_{ms}h_{rj}-h_{mj}h_{rs}+Kg_{ms}g_{rj}-Kg_{mj}g_{rs})h_{i}^{m}\\
		&= F^{rs}(\n^{2}_{rs}h_{ij}+ h_{ms}h^{m}_{r}h_{ij}) +KFg_{ij}-Fh_{mj}h^{m}_{i}-KF^{rs}g_{rs}h_{ij}.
		}
(ii) As in the proof of \autoref{ev-p-parabolic}, let $\eta$ be a smooth lower support of $p$ at $\xi_{0}\in M$ and define $\vp = \eta F$. Then, in coordinates with
\eq{\label{coord} g_{ij}=\de_{ij}, \q h_{ij} = \ka_{i}\de_{ij}}
we have inequality \eqref{Lem5}, which combined with (i) leads to
\eq{-F F^{kl}\n^{2}_{kl}\eta &= -F^{kl}\n^{2}_{kl}\vp+\fr{\vp}{F}F^{kl}\n^{2}_{kl}F+2F^{kl}\n_{k}\eta\n_{l} F\\\
				&\geq 2F^{kl}\n_{k}\eta\n_{l} F + 2\sum_{j>D}\fr{F^{kk}(\n_{k}h^{1}_{j})^{2}}{\ka_{j}-\ka_{1}}+ F^{pq,rs}\n_{1}h_{pq}\n_{1}h_{rs}\\
				&\hp{=}+\eta F^{kl}\n_{kl}^2F- \n_{11}^2 F +F^{ii}(\ka_{i}\ka_{1}+K)(\ka_{i}-\ka_{1}).
	}	
}

In order to prove \autoref{thm:elliptic}, let $F = \ga {\rm c}_{_{K}}^{\al}$, then taking derivatives and using \eqref{Hess-c} we get
\eq{\n^2_{kl}F &= \al\ga {\rm c}_{_{K}}^{\al-1}\n^2_{kl}{\rm c}_{_{K}} + \al(\al-1)\ga {\rm c}_{_{K}}^{\al-2}\n_{k}{\rm c}_{_{K}}\n_{l}{\rm c}_{_{K}}\\
		&= K\Big(\frac{u}{\co} h_{kl}-g_{kl}\Big) \alpha F  + \fr{2 \ep}{F}\n_{k}F \n_{l}F,}
where \eq{\label{def-eps} \ep = \fr{\al-1}{2\al} = \fr{1}{2} - \fr{1}{2\al}.}
This implies
\eq{\label{pf:elliptic-1}p F^{kl}\n^{2}_{kl}F - & \n_{11}^{2}F  = p \alpha K F \Big(\frac{u}{\co} F- F^{kl}g_{kl}\Big) - \alpha K F \Big(\frac{u}{\co} \kappa_1-1\Big)\\
	&+\fr{2 \ep p}{F}F^{kl}\n_{k}F\n_{l}F-\fr{2 \ep}{F}(\n_{1}F)^{2}\\
	=~& K\al F^{ii}(\kappa_i- \kappa_1 )+ \fr{2 \ep}{F}\big(p F^{kl}\n_{k}F\n_{l}F-(\n_{1}F)^{2}\big).}
Here we have used Euler's relation and computations are done in the coordinates \eqref{coord}.

\begin{proof}[Proof of \autoref{thm:elliptic}]
As $|\alpha| \geq 1$, we have that $\ep \in [0,1]$ for $\ep$ defined as in \eqref{def-eps}. Then using that $M$ is convex, we can estimate
\eq{ &2\sum_{j>D}\fr{F^{kk}(\n_{k}h^{1}_{j})^{2}}{\ka_{j}-\ka_{1}} 
			\geq ~&2\ep\sum_{j=1}^n\fr{F^{kk}(\n_{k}h^{1}_{j})^{2}}{\ka_{j}}
			 -2\ep \fr{p}{F}F^{kl}\n_{k}F\n_{l}F + T\ast \nabla p, 
		}
			where we also used the following (see \cite[Lemma~5]{BrendleChoiDaskalopoulos:/2017}):
\eq{\n_{k}\vp\de_{ij} = \n_{k}h_{ij}\q \text{for all} \q 1\leq i,j\leq D.}

Plugging this in \autoref{lem-el} (ii) and using the convexity of $F$, we get 
\eq{- F F^{kl}\n^{2}_{kl}p &\geq  2\ep\sum_{j=1}^n\fr{F^{kk}(\n_{k}h^{1}_{j})^{2}}{\ka_{j}}
	-2\ep \fr{p}{F}F^{kl}\n_{k}F\n_{l}F + \ep F^{pq,rs}\n_{1}h_{pq}\n_{1}h_{rs}\\
	  &\hp{=}  + p F^{kl}\n_{kl}F -\n_{11}F
	 +F^{ii} {\rm Sec}_{i1} (\ka_{i}-\ka_{1}) + T\ast \nabla p\\
	  &\underset{\eqref{pf:elliptic-1}}{\geq} 2\ep\sum_{j=1}^n\fr{F^{kk}(\n_{k}h^{1}_{j})^{2}}{\ka_{j}}
	 + \ep F^{pq,rs}\n_{1}h_{pq}\n_{1}h_{rs} - \frac{2\ep}{F} (\nabla_1 F)^2\\
	  &\hp{=} +F^{ii} ({\rm Sec}_{i1} + \al K) (\ka_{i}-\ka_{1}) + T\ast \nabla p.}
Notice that in the case $\al = -1$, then $\ep = 1$ and we do not need the convexity of $F$ in the first inequality above.

By \autoref{invF} (b) our $F$ is inverse concave, and hence \eqref{inv-conc} leads to
\eq{- F F^{kl}\n^{2}_{kl}p + T \ast \nabla p\geq  F^{ii} ({\rm Sec}_{i1} + \al K) (\ka_{i}-\ka_{1}).}
This completes the proof using the strong maximum principle for viscosity solutions and $M$ has to be centred at the origin provided $K\neq 0$, since only in this case can ${\rm c}_{_{K}}$ be constant.
\end{proof}

\begin{rem}
Note that this approach also provides a direct maximum principle proof of Liebmann's soap bubble theorem (the convex case of Alexandrov's theorem), \cite{Liebmann:03/1900}.
\end{rem}

\subsection*{Solitons}

We complete this paper by proving \autoref{cor:elliptic}.

\begin{proof}[Proof of \autoref{cor:elliptic}]For the given hypersurface $M$, there is a dual hypersurface 
$\ti M\sub \bbM_{K}^{n+1}$, where $K = {\rm sgn }(\mathbb M)$ as in \eqref{signat},  and  
with the properties
\eq{\ti\ka_{i} = \fr{1}{\ka_{i}}, \q \ti {\rm c}_{_{K}} = u,}
see \cite[Thm.~10.4.4, Thm.~10.4.9]{Gerhardt:/2006} and \cite{Scheuer:10/2021}.
Hence $\ti M$ satisfies the equation
\eq{\ti c_{_{K}}^{\fr{1}{\be}} = F(\ka_{i}) = F(\ti\ka_{i}^{-1}) = \fr{1}{\ti F(\ti\ka_{i})},}
i.e. with $\al = -1/\be$ we have
\eq{\ti F_{|\ti M} = \ti {\rm c}_{_{K}}^{\al},}
where $\ti F$ is the inverse curvature function of $F$. Therefore to complete the proof it only remains to check the conditions of \autoref{thm:elliptic} for $\ti M$, which holds provided that for any $\ti g$-orthonormal frame we have
\eq{\widetilde {\rm Sec}_{ij}\geq -\al K.}
In coordinates that diagonalize $\tilde A$, the Gauss equation \eqref{Gauss} for $\ti M$ gives
\eq{\ti R_{ijij} +\al K &= \ti h_{ii}\ti h_{jj} - \ti h_{ij}\ti h_{ij} + (1+\al)K  = \fr{1}{\ka_{i}\ka_{j}} + (1+\al)K
			\geq 0,}
provided
\eq{\fr{1}{\ka_{i}\ka_{j}}\geq \fr{1-\be}{\be}K.}
Notice that if $\fr{1-\be}{\be}K \leq 0$ this is guaranteed by convexity of $M$; otherwise, the inequality follows by the assumption on ${\rm Sec}_M$. Hence the statement follows by direct application of \autoref{thm:elliptic}.
\end{proof}


\providecommand{\bysame}{\leavevmode\hbox to3em{\hrulefill}\thinspace}
\providecommand{\href}[2]{#2}

\end{document}